\documentclass{article}
\long\def\ig#1{\relax}
\ig{Thanks to Roberto Minio for this def'n.  Compare the def'n of
\comment in AMSTeX.}

\newcount \coefa
\newcount \coefb
\newcount \coefc
\newcount\tempcounta
\newcount\tempcountb
\newcount\tempcountc
\newcount\tempcountd
\newcount\xext
\newcount\yext
\newcount\xoff
\newcount\yoff
\newcount\gap%
\newcount\arrowtypea
\newcount\arrowtypeb
\newcount\arrowtypec
\newcount\arrowtyped
\newcount\arrowtypee
\newcount\height
\newcount\width
\newcount\xpos
\newcount\ypos
\newcount\run
\newcount\rise
\newcount\arrowlength
\newcount\halflength
\newcount\arrowtype
\newdimen\tempdimen
\newdimen\xlen
\newdimen\ylen
\newsavebox{\tempboxa}%
\newsavebox{\tempboxb}%
\newsavebox{\tempboxc}%

\makeatletter
\def\settypes(#1,#2,#3){\arrowtypea#1 \arrowtypeb#2 \arrowtypec#3}
\def\settoheight#1#2{\setbox\@tempboxa\hbox{#2}#1\ht\@tempboxa\relax}%
\def\settodepth#1#2{\setbox\@tempboxa\hbox{#2}#1\dp\@tempboxa\relax}%
\def\settokens[#1`#2`#3`#4]{%
     \def\tokena{#1}\def\tokenb{#2}\def\tokenc{#3}\def\tokend{#4}}
\def\setsqparms[#1`#2`#3`#4;#5`#6]{%
\arrowtypea #1
\arrowtypeb #2
\arrowtypec #3
\arrowtyped #4
\width #5
\height #6
}
\def\setpos(#1,#2){\xpos=#1 \ypos#2}

\def\bfig{\begin{picture}(\xext,\yext)(\xoff,\yoff)}
\def\efig{\end{picture}}

\def\putbox(#1,#2)#3{\put(#1,#2){\makebox(0,0){$#3$}}}

\def\settriparms[#1`#2`#3;#4]{\settripairparms[#1`#2`#3`1`1;#4]}%

\def\settripairparms[#1`#2`#3`#4`#5;#6]{%
\arrowtypea #1
\arrowtypeb #2
\arrowtypec #3
\arrowtyped #4
\arrowtypee #5
\width #6
\height #6
}

\def\resetparms{\settripairparms[1`1`1`1`1;500]\width 500}%default values%

\resetparms

\def\mvector(#1,#2)#3{%%
\put(0,0){\vector(#1,#2){#3}}%
\put(0,0){\vector(#1,#2){30}}%
}
\def\evector(#1,#2)#3{{%%
\arrowlength #3
\put(0,0){\vector(#1,#2){\arrowlength}}%
\advance \arrowlength by-30
\put(0,0){\vector(#1,#2){\arrowlength}}%
}}

\def\horsize#1#2{%
\settowidth{\tempdimen}{$#2$}%
#1=\tempdimen
\divide #1 by\unitlength
}

\def\vertsize#1#2{%
\settoheight{\tempdimen}{$#2$}%
#1=\tempdimen
\settodepth{\tempdimen}{$#2$}%
\advance #1 by\tempdimen
\divide #1 by\unitlength
}

\def\vertadjust[#1`#2`#3]{%
\vertsize{\tempcounta}{#1}%
\vertsize{\tempcountb}{#2}%
\ifnum \tempcounta<\tempcountb \tempcounta=\tempcountb \fi
\divide\tempcounta by2
\vertsize{\tempcountb}{#3}%
\ifnum \tempcountb>0 \advance \tempcountb by20 \fi
\ifnum \tempcounta<\tempcountb \tempcounta=\tempcountb \fi
}

\def\horadjust[#1`#2`#3]{%
\horsize{\tempcounta}{#1}%
\horsize{\tempcountb}{#2}%
\ifnum \tempcounta<\tempcountb \tempcounta=\tempcountb \fi
\divide\tempcounta by20
\horsize{\tempcountb}{#3}%
\ifnum \tempcountb>0 \advance \tempcountb by60 \fi
\ifnum \tempcounta<\tempcountb \tempcounta=\tempcountb \fi
}

\ig{ In this procedure, #1 is the paramater that sticks out all the way,
#2 sticks out the least and #3 is a label sticking out half way.  #4 is
the amount of the offset.}

\def\sladjust[#1`#2`#3]#4{%
\tempcountc=#4
\horsize{\tempcounta}{#1}%
\divide \tempcounta by2
\horsize{\tempcountb}{#2}%
\divide \tempcountb by2
\advance \tempcountb by-\tempcountc
\ifnum \tempcounta<\tempcountb \tempcounta=\tempcountb\fi
\divide \tempcountc by2
\horsize{\tempcountb}{#3}%
\advance \tempcountb by-\tempcountc
\ifnum \tempcountb>0 \advance \tempcountb by80\fi
\ifnum \tempcounta<\tempcountb \tempcounta=\tempcountb\fi
\advance\tempcounta by20
}

\def\putvector(#1,#2)(#3,#4)#5#6{{%
\xpos=#1
\ypos=#2
\run=#3
\rise=#4
\arrowlength=#5
\arrowtype=#6
\ifnum \arrowtype<0
    \ifnum \run=0
        \advance \ypos by-\arrowlength
    \else
        \tempcounta \arrowlength
        \multiply \tempcounta by\rise
        \divide \tempcounta by\run
        \ifnum\run>0
            \advance \xpos by\arrowlength
            \advance \ypos by\tempcounta
        \else
            \advance \xpos by-\arrowlength
            \advance \ypos by-\tempcounta
        \fi
    \fi
    \multiply \arrowtype by-1
    \multiply \rise by-1
    \multiply \run by-1
\fi
\ifnum \arrowtype=1
    \put(\xpos,\ypos){\vector(\run,\rise){\arrowlength}}%
\else\ifnum \arrowtype=2
    \put(\xpos,\ypos){\mvector(\run,\rise)\arrowlength}%
\else\ifnum\arrowtype=3
    \put(\xpos,\ypos){\evector(\run,\rise){\arrowlength}}%
\fi\fi\fi
}}

\def\putsplitvector(#1,#2)#3#4{%%
\xpos #1
\ypos #2
\arrowtype #4
\halflength #3
\arrowlength #3
\gap 140
\advance \halflength by-\gap
\divide \halflength by2
\ifnum \arrowtype=1
    \put(\xpos,\ypos){\line(0,-1){\halflength}}%
    \advance\ypos by-\halflength
    \advance\ypos by-\gap
    \put(\xpos,\ypos){\vector(0,-1){\halflength}}%
\else\ifnum \arrowtype=2
    \put(\xpos,\ypos){\line(0,-1)\halflength}%
    \put(\xpos,\ypos){\vector(0,-1)3}%
    \advance\ypos by-\halflength
    \advance\ypos by-\gap
    \put(\xpos,\ypos){\vector(0,-1){\halflength}}%
\else\ifnum\arrowtype=3
    \put(\xpos,\ypos){\line(0,-1)\halflength}%
    \advance\ypos by-\halflength
    \advance\ypos by-\gap
    \put(\xpos,\ypos){\evector(0,-1){\halflength}}%
\else\ifnum \arrowtype=-1
    \advance \ypos by-\arrowlength
    \put(\xpos,\ypos){\line(0,1){\halflength}}%
    \advance\ypos by\halflength
    \advance\ypos by\gap
    \put(\xpos,\ypos){\vector(0,1){\halflength}}%
\else\ifnum \arrowtype=-2
    \advance \ypos by-\arrowlength
    \put(\xpos,\ypos){\line(0,1)\halflength}%
    \put(\xpos,\ypos){\vector(0,1)3}%
    \advance\ypos by\halflength
    \advance\ypos by\gap
    \put(\xpos,\ypos){\vector(0,1){\halflength}}%
\else\ifnum\arrowtype=-3
    \advance \ypos by-\arrowlength
    \put(\xpos,\ypos){\line(0,1)\halflength}%
    \advance\ypos by\halflength
    \advance\ypos by\gap
    \put(\xpos,\ypos){\evector(0,1){\halflength}}%
\fi\fi\fi\fi\fi\fi
}

\def\putmorphism(#1)(#2,#3)[#4`#5`#6]#7#8#9{{%
\run #2
\rise #3
\ifnum\rise=0
  \puthmorphism(#1)[#4`#5`#6]{#7}{#8}{#9}%
\else\ifnum\run=0
  \putvmorphism(#1)[#4`#5`#6]{#7}{#8}{#9}%
\else
\setpos(#1)%
\arrowlength #7
\arrowtype #8
\ifnum\run=0
\else\ifnum\rise=0
\else
\ifnum\run>0
    \coefa=1
\else
   \coefa=-1
\fi
\ifnum\arrowtype>0
   \coefb=0
   \coefc=-1
\else
   \coefb=\coefa
   \coefc=1
   \arrowtype=-\arrowtype
\fi
\width=2
\multiply \width by\run
\divide \width by\rise
\ifnum \width<0  \width=-\width\fi
\advance\width by60
\if l#9 \width=-\width\fi
\putbox(\xpos,\ypos){#4}%            %node 1
{\multiply \coefa by\arrowlength%      %node 2
\advance\xpos by\coefa
\multiply \coefa by\rise
\divide \coefa by\run
\advance \ypos by\coefa
\putbox(\xpos,\ypos){#5} }%
{\multiply \coefa by\arrowlength%      %label
\divide \coefa by2
\advance \xpos by\coefa
\advance \xpos by\width
\multiply \coefa by\rise
\divide \coefa by\run
\advance \ypos by\coefa
\if l#9%
   \put(\xpos,\ypos){\makebox(0,0)[r]{$#6$}}%
\else\if r#9%
   \put(\xpos,\ypos){\makebox(0,0)[l]{$#6$}}%
\fi\fi }%
{\multiply \rise by-\coefc%             %arrow
\multiply \run by-\coefc
\multiply \coefb by\arrowlength
\advance \xpos by\coefb
\multiply \coefb by\rise
\divide \coefb by\run
\advance \ypos by\coefb
\multiply \coefc by70
\advance \ypos by\coefc
\multiply \coefc by\run
\divide \coefc by\rise
\advance \xpos by\coefc
\multiply \coefa by140
\multiply \coefa by\run
\divide \coefa by\rise
\advance \arrowlength by\coefa
\ifnum \arrowtype=1
   \put(\xpos,\ypos){\vector(\run,\rise){\arrowlength}}%
\else\ifnum\arrowtype=2
   \put(\xpos,\ypos){\mvector(\run,\rise){\arrowlength}}%
\else\ifnum\arrowtype=3
   \put(\xpos,\ypos){\evector(\run,\rise){\arrowlength}}%
\fi\fi\fi}\fi\fi\fi\fi}}

\def\puthmorphism(#1,#2)[#3`#4`#5]#6#7#8{{%
\xpos #1
\ypos #2
\width #6
\arrowlength #6
\putbox(\xpos,\ypos){#3\vphantom{#4}}%
{\advance \xpos by\arrowlength
\putbox(\xpos,\ypos){\vphantom{#3}#4}}%
\horsize{\tempcounta}{#3}%
\horsize{\tempcountb}{#4}%
\divide \tempcounta by2
\divide \tempcountb by2
\advance \tempcounta by30
\advance \tempcountb by30
\advance \xpos by\tempcounta
\advance \arrowlength by-\tempcounta
\advance \arrowlength by-\tempcountb
\putvector(\xpos,\ypos)(1,0){\arrowlength}{#7}%
\divide \arrowlength by2
\advance \xpos by\arrowlength
\vertsize{\tempcounta}{#5}%
\divide\tempcounta by2
\advance \tempcounta by20
\if a#8 %
   \advance \ypos by\tempcounta
   \putbox(\xpos,\ypos){#5}%
\else
   \advance \ypos by-\tempcounta
   \putbox(\xpos,\ypos){#5}%
\fi}}

\def\putvmorphism(#1,#2)[#3`#4`#5]#6#7#8{{%
\xpos #1
\ypos #2
\arrowlength #6
\arrowtype #7
\settowidth{\xlen}{$#5$}%
\putbox(\xpos,\ypos){#3}%
{\advance \ypos by-\arrowlength
\putbox(\xpos,\ypos){#4}}%
{\advance\arrowlength by-140
\advance \ypos by-70
\ifdim\xlen>0pt
   \if m#8%
      \putsplitvector(\xpos,\ypos){\arrowlength}{\arrowtype}%
   \else
      \putvector(\xpos,\ypos)(0,-1){\arrowlength}{\arrowtype}%
   \fi
\else
   \putvector(\xpos,\ypos)(0,-1){\arrowlength}{\arrowtype}%
\fi}%
\ifdim\xlen>0pt
   \divide \arrowlength by2
   \advance\ypos by-\arrowlength
   \if l#8%
      \advance \xpos by-40
      \put(\xpos,\ypos){\makebox(0,0)[r]{$#5$}}%
   \else\if r#8%
      \advance \xpos by40
      \put(\xpos,\ypos){\makebox(0,0)[l]{$#5$}}%
   \else
      \putbox(\xpos,\ypos){#5}%
   \fi\fi
\fi
}}

\def\topadjust[#1`#2`#3]{%
\yoff=10
\vertadjust[#1`#2`{#3}]%
\advance \yext by\tempcounta
\advance \yext by 10
}
\def\botadjust[#1`#2`#3]{%
\vertadjust[#1`#2`{#3}]%
\advance \yext by\tempcounta
\advance \yoff by-\tempcounta
}
\def\leftadjust[#1`#2`#3]{%
\xoff=0
\horadjust[#1`#2`{#3}]%
\advance \xext by\tempcounta
\advance \xoff by-\tempcounta
}
\def\rightadjust[#1`#2`#3]{%
\horadjust[#1`#2`{#3}]%
\advance \xext by\tempcounta
}
\def\rightsladjust[#1`#2`#3]{%
\sladjust[#1`#2`{#3}]{\width}%
\advance \xext by\tempcounta
}
\def\leftsladjust[#1`#2`#3]{%
\xoff=0
\sladjust[#1`#2`{#3}]{\width}%
\advance \xext by\tempcounta
\advance \xoff by-\tempcounta
}
\def\adjust[#1`#2;#3`#4;#5`#6;#7`#8]{%
\topadjust[#1``{#2}]
\leftadjust[#3``{#4}]
\rightadjust[#5``{#6}]
\botadjust[#7``{#8}]}

\def\putsquarep<#1>(#2)[#3;#4`#5`#6`#7]{{%
\setsqparms[#1]%
\setpos(#2)%
\settokens[#3]%
\puthmorphism(\xpos,\ypos)[\tokenc`\tokend`{#7}]{\width}{\arrowtyped}b%
\advance\ypos by \height
\puthmorphism(\xpos,\ypos)[\tokena`\tokenb`{#4}]{\width}{\arrowtypea}a%
\putvmorphism(\xpos,\ypos)[``{#5}]{\height}{\arrowtypeb}l%
\advance\xpos by \width
\putvmorphism(\xpos,\ypos)[``{#6}]{\height}{\arrowtypec}r%
}}

\def\putsquare{\@ifnextchar <{\putsquarep}{\putsquarep%
   <\arrowtypea`\arrowtypeb`\arrowtypec`\arrowtyped;\width`\height>}}
\def\square{\@ifnextchar< {\squarep}{\squarep
   <\arrowtypea`\arrowtypeb`\arrowtypec`\arrowtyped;\width`\height>}}
\def\squarep<#1>[#2`#3`#4`#5;#6`#7`#8`#9]{{%          %     #2------>#3
\setsqparms[#1]%                                      %      |       |
\xext=\width                                          %      |       |
\yext=\height                                         %    #7|       |#8
\topadjust[#2`#3`{#6}]%                               %      |       |
\botadjust[#4`#5`{#9}]%                               %      |       |
\leftadjust[#2`#4`{#7}]%                              %
\rightadjust[#3`#5`{#8}]%                             %     #4------>#5
\begin{picture}(\xext,\yext)(\xoff,\yoff)%                      #9
\putsquarep<\arrowtypea`\arrowtypeb`\arrowtypec`\arrowtyped;\width`\height>%
(0,0)[#2`#3`#4`#5;#6`#7`#8`{#9}]%
\end{picture}%
}}

\def\putptrianglep<#1>(#2,#3)[#4`#5`#6;#7`#8`#9]{{%
\settriparms[#1]%
\xpos=#2 \ypos=#3
\advance\ypos by \height
\puthmorphism(\xpos,\ypos)[#4`#5`{#7}]{\height}{\arrowtypea}a%
\putvmorphism(\xpos,\ypos)[`#6`{#8}]{\height}{\arrowtypeb}l%
\advance\xpos by\height
\putmorphism(\xpos,\ypos)(-1,-1)[``{#9}]{\height}{\arrowtypec}r%
}}

\def\putptriangle{\@ifnextchar <{\putptrianglep}{\putptrianglep
   <\arrowtypea`\arrowtypeb`\arrowtypec;\height>}}
\def\ptriangle{\@ifnextchar <{\ptrianglep}{\ptrianglep
   <\arrowtypea`\arrowtypeb`\arrowtypec;\height>}}

\def\ptrianglep<#1>[#2`#3`#4;#5`#6`#7]{{%%       #5
\settriparms[#1]%
\width=\height                         %      #2----->#3
\xext=\width                           %      |      /
\yext=\width                           %      |     /
\topadjust[#2`#3`{#5}]%                %    #6|    /#7
\botadjust[#3``]%                      %      |   /
\leftadjust[#2`#4`{#6}]%               %      |  /
\rightsladjust[#3`#4`{#7}]%            %
\begin{picture}(\xext,\yext)(\xoff,\yoff)%    #4
\putptrianglep<\arrowtypea`\arrowtypeb`\arrowtypec;\height>%
(0,0)[#2`#3`#4;#5`#6`{#7}]%
\end{picture}%
}}

\def\putqtrianglep<#1>(#2,#3)[#4`#5`#6;#7`#8`#9]{{%
\settriparms[#1]%
\xpos=#2 \ypos=#3
\advance\ypos by\height
\puthmorphism(\xpos,\ypos)[#4`#5`{#7}]{\height}{\arrowtypea}a%
\putmorphism(\xpos,\ypos)(1,-1)[``{#8}]{\height}{\arrowtypeb}l%
\advance\xpos by\height
\putvmorphism(\xpos,\ypos)[`#6`{#9}]{\height}{\arrowtypec}r%
}}

\def\putqtriangle{\@ifnextchar <{\putqtrianglep}{\putqtrianglep
   <\arrowtypea`\arrowtypeb`\arrowtypec;\height>}}
\def\qtriangle{\@ifnextchar <{\qtrianglep}{\qtrianglep
   <\arrowtypea`\arrowtypeb`\arrowtypec;\height>}}

\def\qtrianglep<#1>[#2`#3`#4;#5`#6`#7]{{%%
\settriparms[#1]%                                  #5
\width=\height                         %        #2----->#3
\xext=\width                           %         \      |
\yext=\height                          %          \     |
\topadjust[#2`#3`{#5}]%                %         #6\    |#7
\botadjust[#4``]%                      %            \   |
\leftsladjust[#2`#4`{#6}]%             %             \  |
\rightadjust[#3`#4`{#7}]%              %
\begin{picture}(\xext,\yext)(\xoff,\yoff)%             #4
\putqtrianglep<\arrowtypea`\arrowtypeb`\arrowtypec;\height>%
(0,0)[#2`#3`#4;#5`#6`{#7}]%
\end{picture}%
}}

\def\putdtrianglep<#1>(#2,#3)[#4`#5`#6;#7`#8`#9]{{%
\settriparms[#1]%
\xpos=#2 \ypos=#3
\puthmorphism(\xpos,\ypos)[#5`#6`{#9}]{\height}{\arrowtypec}b%
\advance\xpos by \height \advance\ypos by\height
\putmorphism(\xpos,\ypos)(-1,-1)[``{#7}]{\height}{\arrowtypea}l%
\putvmorphism(\xpos,\ypos)[#4``{#8}]{\height}{\arrowtypeb}r%
}}

\def\putdtriangle{\@ifnextchar <{\putdtrianglep}{\putdtrianglep
   <\arrowtypea`\arrowtypeb`\arrowtypec;\height>}}
\def\dtriangle{\@ifnextchar <{\dtrianglep}{\dtrianglep
   <\arrowtypea`\arrowtypeb`\arrowtypec;\height>}}

\def\dtrianglep<#1>[#2`#3`#4;#5`#6`#7]{{%%
\settriparms[#1]%                                          #2
\width=\height                         %                  / |
\xext=\width                           %                 /  |
\yext=\height                          %              #5/   |#6
\topadjust[#2``]%                      %               /    |
\botadjust[#3`#4`{#7}]%                %              /     |
\leftsladjust[#3`#2`{#5}]%             %
\rightadjust[#2`#4`{#6}]%              %            #3----->#4
\begin{picture}(\xext,\yext)(\xoff,\yoff)%              #7
\putdtrianglep<\arrowtypea`\arrowtypeb`\arrowtypec;\height>%
(0,0)[#2`#3`#4;#5`#6`{#7}]%
\end{picture}%
}}

\def\putbtrianglep<#1>(#2,#3)[#4`#5`#6;#7`#8`#9]{{%
\settriparms[#1]%
\xpos=#2 \ypos=#3
\puthmorphism(\xpos,\ypos)[#5`#6`{#9}]{\height}{\arrowtypec}b%
\advance\ypos by\height
\putmorphism(\xpos,\ypos)(1,-1)[``{#8}]{\height}{\arrowtypeb}r%
\putvmorphism(\xpos,\ypos)[#4``{#7}]{\height}{\arrowtypea}l%
}}

\def\putbtriangle{\@ifnextchar <{\putbtrianglep}{\putbtrianglep
   <\arrowtypea`\arrowtypeb`\arrowtypec;\height>}}
\def\btriangle{\@ifnextchar <{\btrianglep}{\btrianglep
   <\arrowtypea`\arrowtypeb`\arrowtypec;\height>}}

\def\btrianglep<#1>[#2`#3`#4;#5`#6`#7]{{%%
\settriparms[#1]%                                     #2
\width=\height                         %              | \
\xext=\width                           %              |  \
\yext=\height                          %            #5|   \#6
\topadjust[#2``]%                      %              |    \
\botadjust[#3`#4`{#7}]%                %              |     \
\leftadjust[#2`#3`{#5}]%               %
\rightsladjust[#4`#2`{#6}]%            %              #3----->#4
\begin{picture}(\xext,\yext)(\xoff,\yoff)%                #7
\putbtrianglep<\arrowtypea`\arrowtypeb`\arrowtypec;\height>%
(0,0)[#2`#3`#4;#5`#6`{#7}]%
\end{picture}%
}}

\def\putAtrianglep<#1>(#2,#3)[#4`#5`#6;#7`#8`#9]{{%
\settriparms[#1]%
\xpos=#2 \ypos=#3
{\multiply \height by2
\puthmorphism(\xpos,\ypos)[#5`#6`{#9}]{\height}{\arrowtypec}b}%
\advance\xpos by\height \advance\ypos by\height
\putmorphism(\xpos,\ypos)(-1,-1)[#4``{#7}]{\height}{\arrowtypea}l%
\putmorphism(\xpos,\ypos)(1,-1)[``{#8}]{\height}{\arrowtypeb}r%
}}

\def\putAtriangle{\@ifnextchar <{\putAtrianglep}{\putAtrianglep
   <\arrowtypea`\arrowtypeb`\arrowtypec;\height>}}
\def\Atriangle{\@ifnextchar <{\Atrianglep}{\Atrianglep
   <\arrowtypea`\arrowtypeb`\arrowtypec;\height>}}

\def\Atrianglep<#1>[#2`#3`#4;#5`#6`#7]{{%%
\settriparms[#1]%                                 #2
\width=\height                         %         /   \
\xext=\width                           %        /     \
\yext=\height                          %     #5/       \#6
\topadjust[#2``]%                      %      /         \
\botadjust[#3`#4`{#7}]%                %     /           \
\multiply \xext by2 %                  %
\leftsladjust[#3`#2`{#5}]%             %   #3------------>#4
\rightsladjust[#4`#2`{#6}]%            %          #7
\begin{picture}(\xext,\yext)(\xoff,\yoff)%
\putAtrianglep<\arrowtypea`\arrowtypeb`\arrowtypec;\height>%
(0,0)[#2`#3`#4;#5`#6`{#7}]%
\end{picture}%
}}

\def\putAtrianglepairp<#1>(#2)[#3;#4`#5`#6`#7`#8]{{
\settripairparms[#1]%
\setpos(#2)%
\settokens[#3]%
\puthmorphism(\xpos,\ypos)[\tokenb`\tokenc`{#7}]{\height}{\arrowtyped}b%
\advance\xpos by\height
\advance\ypos by\height
\putmorphism(\xpos,\ypos)(-1,-1)[\tokena``{#4}]{\height}{\arrowtypea}l%
\putvmorphism(\xpos,\ypos)[``{#5}]{\height}{\arrowtypeb}m%
\putmorphism(\xpos,\ypos)(1,-1)[``{#6}]{\height}{\arrowtypec}r%
}}

\def\putAtrianglepair{\@ifnextchar <{\putAtrianglepairp}{\putAtrianglepairp%
   <\arrowtypea`\arrowtypeb`\arrowtypec`\arrowtyped`\arrowtypee;\height>}}
\def\Atrianglepair{\@ifnextchar <{\Atrianglepairp}{\Atrianglepairp%
   <\arrowtypea`\arrowtypeb`\arrowtypec`\arrowtyped`\arrowtypee;\height>}}

\def\Atrianglepairp<#1>[#2;#3`#4`#5`#6`#7]{{%
\settripairparms[#1]%
\settokens[#2]%
\width=\height
\xext=\width
\yext=\height
\topadjust[\tokena``]%
\vertadjust[\tokenb`\tokenc`{#6}]%                      %  #2a
\tempcountd=\tempcounta                       %           / | \
\vertadjust[\tokenc`\tokend`{#7}]%            %          /  |  \
\ifnum\tempcounta<\tempcountd                 %       #3/  #4   \#5
\tempcounta=\tempcountd\fi                    %        /    |    \
\advance \yext by\tempcounta                  %       /     |     \
\advance \yoff by-\tempcounta                 %
\multiply \xext by2 %                         %     #2b---->#2c---->#2d
\leftsladjust[\tokenb`\tokena`{#3}]%          %         #6     #7
\rightsladjust[\tokend`\tokena`{#5}]%
\begin{picture}(\xext,\yext)(\xoff,\yoff)%
\putAtrianglepairp
<\arrowtypea`\arrowtypeb`\arrowtypec`\arrowtyped`\arrowtypee;\height>%
(0,0)[#2;#3`#4`#5`#6`{#7}]%
\end{picture}%
}}

\def\putVtrianglep<#1>(#2,#3)[#4`#5`#6;#7`#8`#9]{{%
\settriparms[#1]%
\xpos=#2 \ypos=#3
\advance\ypos by\height
{\multiply\height by2
\puthmorphism(\xpos,\ypos)[#4`#5`{#7}]{\height}{\arrowtypea}a}%
\putmorphism(\xpos,\ypos)(1,-1)[`#6`{#8}]{\height}{\arrowtypeb}l%
\advance\xpos by\height
\advance\xpos by\height
\putmorphism(\xpos,\ypos)(-1,-1)[``{#9}]{\height}{\arrowtypec}r%
}}

\def\putVtriangle{\@ifnextchar <{\putVtrianglep}{\putVtrianglep
   <\arrowtypea`\arrowtypeb`\arrowtypec;\height>}}
\def\Vtriangle{\@ifnextchar <{\Vtrianglep}{\Vtrianglep
   <\arrowtypea`\arrowtypeb`\arrowtypec;\height>}}

\def\Vtrianglep<#1>[#2`#3`#4;#5`#6`#7]{{%%
\settriparms[#1]%                                      #5
\width=\height                         %        #2------------->#3
\xext=\width                           %         \             /
\yext=\height                          %          \           /
\topadjust[#2`#3`{#5}]%                %         #6\         /#7
\botadjust[#4``]%                      %            \       /
\multiply \xext by2 %                  %             \     /
\leftsladjust[#2`#3`{#6}]%             %
\rightsladjust[#3`#4`{#7}]%            %               #4
\begin{picture}(\xext,\yext)(\xoff,\yoff)%
\putVtrianglep<\arrowtypea`\arrowtypeb`\arrowtypec;\height>%
(0,0)[#2`#3`#4;#5`#6`{#7}]%
\end{picture}%
}}

\def\putVtrianglepairp<#1>(#2)[#3;#4`#5`#6`#7`#8]{{
\settripairparms[#1]%
\setpos(#2)%
\settokens[#3]%
\advance\ypos by\height
\putmorphism(\xpos,\ypos)(1,-1)[`\tokend`{#6}]{\height}{\arrowtypec}l%
\puthmorphism(\xpos,\ypos)[\tokena`\tokenb`{#4}]{\height}{\arrowtypea}a%
\advance\xpos by\height
\putvmorphism(\xpos,\ypos)[``{#7}]{\height}{\arrowtyped}m%
\advance\xpos by\height
\putmorphism(\xpos,\ypos)(-1,-1)[``{#8}]{\height}{\arrowtypee}r%
}}

\def\putVtrianglepair{\@ifnextchar <{\putVtrianglepairp}{\putVtrianglepairp%
    <\arrowtypea`\arrowtypeb`\arrowtypec`\arrowtyped`\arrowtypee;\height>}}
\def\Vtrianglepair{\@ifnextchar <{\Vtrianglepairp}{\Vtrianglepairp%
    <\arrowtypea`\arrowtypeb`\arrowtypec`\arrowtyped`\arrowtypee;\height>}}

\def\Vtrianglepairp<#1>[#2;#3`#4`#5`#6`#7]{{%
\settripairparms[#1]%
\settokens[#2]%                            #3      #4
\xext=\height                  %        #2a---->#2b---->#2c
\width=\height                 %         \      |      /
\yext=\height                  %          \     |     /
\vertadjust[\tokena`\tokenb`{#4}]%       #5\   #6    /#7
\tempcountd=\tempcounta        %            \   |   /
\vertadjust[\tokenb`\tokenc`{#5}]%           \  |  /
\ifnum\tempcounta<\tempcountd%
\tempcounta=\tempcountd\fi%                    #2d
\advance \yext by\tempcounta
\botadjust[\tokend``]%
\multiply \xext by2
\leftsladjust[\tokena`\tokend`{#6}]%
\rightsladjust[\tokenc`\tokend`{#7}]%
\begin{picture}(\xext,\yext)(\xoff,\yoff)%
\putVtrianglepairp
<\arrowtypea`\arrowtypeb`\arrowtypec`\arrowtyped`\arrowtypee;\height>%
(0,0)[#2;#3`#4`#5`#6`{#7}]%
\end{picture}%
}}

\def\putCtrianglep<#1>(#2,#3)[#4`#5`#6;#7`#8`#9]{{%
\settriparms[#1]%
\xpos=#2 \ypos=#3
\advance\ypos by\height
\putmorphism(\xpos,\ypos)(1,-1)[``{#9}]{\height}{\arrowtypec}l%
\advance\xpos by\height
\advance\ypos by\height
\putmorphism(\xpos,\ypos)(-1,-1)[#4`#5`{#7}]{\height}{\arrowtypea}l%
{\multiply\height by 2
\putvmorphism(\xpos,\ypos)[`#6`{#8}]{\height}{\arrowtypeb}r}%
}}

\def\putCtriangle{\@ifnextchar <{\putCtrianglep}{\putCtrianglep
    <\arrowtypea`\arrowtypeb`\arrowtypec;\height>}}
\def\Ctriangle{\@ifnextchar <{\Ctrianglep}{\Ctrianglep
    <\arrowtypea`\arrowtypeb`\arrowtypec;\height>}}

\def\Ctrianglep<#1>[#2`#3`#4;#5`#6`#7]{{%%
\settriparms[#1]%                                         #2
\width=\height                          %                / |
\xext=\width                            %               /  |
\yext=\height                           %            #5/   |
\multiply \yext by2 %                   %             /    |
\topadjust[#2``]%                       %            /     |
\botadjust[#4``]%                       %           v      |
\sladjust[#3`#2`{#5}]{\width}%          %          #3      |#6
\tempcountd=\tempcounta                 %           \      |
\sladjust[#3`#4`{#7}]{\width}%          %            \     |
\ifnum \tempcounta<\tempcountd          %           #7\    |
\tempcounta=\tempcountd\fi              %              \   |
\advance \xext by\tempcounta            %               \  |
\advance \xoff by-\tempcounta           %
\rightadjust[#2`#4`{#6}]%               %                 #4
\begin{picture}(\xext,\yext)(\xoff,\yoff)%
\putCtrianglep<\arrowtypea`\arrowtypeb`\arrowtypec;\height>%
(0,0)[#2`#3`#4;#5`#6`{#7}]%
\end{picture}%
}}

\def\putDtrianglep<#1>(#2,#3)[#4`#5`#6;#7`#8`#9]{{%
\settriparms[#1]%
\xpos=#2 \ypos=#3
\advance\xpos by\height \advance\ypos by\height
\putmorphism(\xpos,\ypos)(-1,-1)[``{#9}]{\height}{\arrowtypec}r%
\advance\xpos by-\height \advance\ypos by\height
\putmorphism(\xpos,\ypos)(1,-1)[`#5`{#8}]{\height}{\arrowtypeb}r%
{\multiply\height by 2
\putvmorphism(\xpos,\ypos)[#4`#6`{#7}]{\height}{\arrowtypea}l}%
}}

\def\putDtriangle{\@ifnextchar <{\putDtrianglep}{\putDtrianglep
    <\arrowtypea`\arrowtypeb`\arrowtypec;\height>}}
\def\Dtriangle{\@ifnextchar <{\Dtrianglep}{\Dtrianglep
   <\arrowtypea`\arrowtypeb`\arrowtypec;\height>}}

\def\Dtrianglep<#1>[#2`#3`#4;#5`#6`#7]{{%%
\settriparms[#1]%                                 #2
\width=\height                         %          | \
\xext=\height                          %          |  \
\yext=\height                          %          |   \#6
\multiply \yext by2 %                  %          |    \
\topadjust[#2``]%                      %          |     \
\botadjust[#4``]%                      %          |
\leftadjust[#2`#4`{#5}]%               %        #5|      #3
\sladjust[#3`#2`{#5}]{\height}%        %          |      /
\tempcountd=\tempcountd                %          |     /
\sladjust[#3`#4`{#7}]{\height}%        %          |    /#7
\ifnum \tempcounta<\tempcountd         %          |   /
\tempcounta=\tempcountd\fi             %          |  /
\advance \xext by\tempcounta           %
\begin{picture}(\xext,\yext)(\xoff,\yoff)%        #4
\putDtrianglep<\arrowtypea`\arrowtypeb`\arrowtypec;\height>%
(0,0)[#2`#3`#4;#5`#6`{#7}]%
\end{picture}%
}}

\def\setrecparms[#1`#2]{\width=#1 \height=#2}%
\def\recursep<#1`#2>[#3;#4`#5`#6`#7`#8]{{%
\width=#1 \height=#2
\settokens[#3]
\settowidth{\tempdimen}{$\tokena$}
\ifdim\tempdimen=0pt
  \savebox{\tempboxa}{\hbox{$\tokenb$}}%
  \savebox{\tempboxb}{\hbox{$\tokend$}}%
  \savebox{\tempboxc}{\hbox{$#6$}}%
\else
  \savebox{\tempboxa}{\hbox{$\hbox{$\tokena$}\times\hbox{$\tokenb$}$}}%
  \savebox{\tempboxb}{\hbox{$\hbox{$\tokena$}\times\hbox{$\tokend$}$}}%
  \savebox{\tempboxc}{\hbox{$\hbox{$\tokena$}\times\hbox{$#6$}$}}%
\fi
\ypos=\height
\divide\ypos by 2
\xpos=\ypos
\advance\xpos by \width
\xext=\xpos \yext=\height
\topadjust[#3`\usebox{\tempboxa}`{#4}]%
\botadjust[#5`\usebox{\tempboxb}`{#8}]%
\sladjust[\tokenc`\tokenb`{#5}]{\ypos}%
\tempcountd=\tempcounta
\sladjust[\tokenc`\tokend`{#5}]{\ypos}%
\ifnum \tempcounta<\tempcountd
\tempcounta=\tempcountd\fi
\advance \xext by\tempcounta
\advance \xoff by-\tempcounta
\rightadjust[\usebox{\tempboxa}`\usebox{\tempboxb}`\usebox{\tempboxc}]%
\bfig
\putCtrianglep<-1`1`1;\ypos>(0,0)[`\tokenc`;#5`#6`{#7}]%
\puthmorphism(\ypos,0)[\tokend`\usebox{\tempboxb}`{#8}]{\width}{-1}b%
\puthmorphism(\ypos,\height)[\tokenb`\usebox{\tempboxa}`{#4}]{\width}{-1}a%
\advance\ypos by \width
\putvmorphism(\ypos,\height)[``\usebox{\tempboxc}]{\height}1r%
\efig
}}

\def\recurse{\@ifnextchar <{\recursep}{\recursep<\width`\height>}}

\def\puttwohmorphisms(#1,#2)[#3`#4;#5`#6]#7#8#9{{%
\puthmorphism(#1,#2)[#3`#4`]{#7}0a
\ypos=#2
\advance\ypos by 20
\puthmorphism(#1,\ypos)[\phantom{#3}`\phantom{#4}`#5]{#7}{#8}a
\advance\ypos by -40
\puthmorphism(#1,\ypos)[\phantom{#3}`\phantom{#4}`#6]{#7}{#9}b
}}

\def\puttwovmorphisms(#1,#2)[#3`#4;#5`#6]#7#8#9{{%
\putvmorphism(#1,#2)[#3`#4`]{#7}0a
\xpos=#1
\advance\xpos by -20
\putvmorphism(\xpos,#2)[\phantom{#3}`\phantom{#4}`#5]{#7}{#8}l
\advance\xpos by 40
\putvmorphism(\xpos,#2)[\phantom{#3}`\phantom{#4}`#6]{#7}{#9}r
}}

\def\puthcoequalizer(#1)[#2`#3`#4;#5`#6`#7]#8#9{{%
\setpos(#1)%
\puttwohmorphisms(\xpos,\ypos)[#2`#3;#5`#6]{#8}11%
\advance\xpos by #8
\puthmorphism(\xpos,\ypos)[\phantom{#3}`#4`#7]{#8}1{#9}
}}

\def\putvcoequalizer(#1)[#2`#3`#4;#5`#6`#7]#8#9{{%
\setpos(#1)%
\puttwovmorphisms(\xpos,\ypos)[#2`#3;#5`#6]{#8}11%
\advance\ypos by -#8
\putvmorphism(\xpos,\ypos)[\phantom{#3}`#4`#7]{#8}1{#9}
}}

\def\putthreehmorphisms(#1)[#2`#3;#4`#5`#6]#7(#8)#9{{%
\setpos(#1) \settypes(#8)
\if a#9 %
     \vertsize{\tempcounta}{#5}%
     \vertsize{\tempcountb}{#6}%
     \ifnum \tempcounta<\tempcountb \tempcounta=\tempcountb \fi
\else
     \vertsize{\tempcounta}{#4}%
     \vertsize{\tempcountb}{#5}%
     \ifnum \tempcounta<\tempcountb \tempcounta=\tempcountb \fi
\fi
\advance \tempcounta by 60
\puthmorphism(\xpos,\ypos)[#2`#3`#5]{#7}{\arrowtypeb}{#9}
\advance\ypos by \tempcounta
\puthmorphism(\xpos,\ypos)[\phantom{#2}`\phantom{#3}`#4]{#7}{\arrowtypea}{#9}
\advance\ypos by -\tempcounta \advance\ypos by -\tempcounta
\puthmorphism(\xpos,\ypos)[\phantom{#2}`\phantom{#3}`#6]{#7}{\arrowtypec}{#9}
}}

\def\putarc(#1,#2)[#3`#4`#5]#6#7#8{{%
\xpos #1
\ypos #2
\width #6
\arrowlength #6
\putbox(\xpos,\ypos){#3\vphantom{#4}}%
{\advance \xpos by\arrowlength
\putbox(\xpos,\ypos){\vphantom{#3}#4}}%
\horsize{\tempcounta}{#3}%
\horsize{\tempcountb}{#4}%
\divide \tempcounta by2
\divide \tempcountb by2
\advance \tempcounta by30
\advance \tempcountb by30
\advance \xpos by\tempcounta
\advance \arrowlength by-\tempcounta
\advance \arrowlength by-\tempcountb
\halflength=\arrowlength \divide\halflength by 2
\divide\arrowlength by 5
\put(\xpos,\ypos){\bezier{\arrowlength}(0,0)(50,50)(\halflength,50)}
\ifnum #7=-1 \put(\xpos,\ypos){\vector(-3,-2)0} \fi
\advance\xpos by \halflength
\put(\xpos,\ypos){\xpos=\halflength \advance\xpos by -50
   \bezier{\arrowlength}(0,50)(\xpos,50)(\halflength,0)}
\ifnum #7=1 {\advance \xpos by
   \halflength \put(\xpos,\ypos){\vector(3,-2)0}} \fi
\advance\ypos by 50
\vertsize{\tempcounta}{#5}%
\divide\tempcounta by2
\advance \tempcounta by20
\if a#8 %
   \advance \ypos by\tempcounta
   \putbox(\xpos,\ypos){#5}%
\else
   \advance \ypos by-\tempcounta
   \putbox(\xpos,\ypos){#5}%
\fi
}}

\makeatother

\usepackage{amsthm}
\usepackage{dsfont}
\usepackage{stmaryrd}
\usepackage[all]{xy}
\usepackage{tikz}
\usetikzlibrary{matrix,arrows}
\advance\hoffset-2truecm
\advance\voffset-1.5truecm
\newtheorem{theorem}{Theorem}[section]

\begin{document}

\sloppy

%commands
\newcommand{\nl}{\hspace{2cm}\\ }

\def\nec{\Box}
\def\pos{\Diamond}
\def\diam{{\tiny\Diamond}}

\def\lc{\lceil}
\def\rc{\rceil}
\def\lf{\lfloor}
\def\rf{\rfloor}
\def\lk{\langle}
\def\rk{\rangle}
\def\blk{\dot{\langle\!\!\langle}}
\def\brk{\dot{\rangle\!\!\rangle}}

\newcommand{\pa}{\parallel}
\newcommand{\lra}{\longrightarrow}
\newcommand{\hra}{\hookrightarrow}
\newcommand{\hla}{\hookleftarrow}
\newcommand{\ra}{\rightarrow}
\newcommand{\la}{\leftarrow}
\newcommand{\lla}{\longleftarrow}
\newcommand{\da}{\downarrow}
\newcommand{\ua}{\uparrow}
\newcommand{\dA}{\downarrow\!\!\!^\bullet}
\newcommand{\uA}{\uparrow\!\!\!_\bullet}
\newcommand{\Da}{\Downarrow}
\newcommand{\DA}{\Downarrow\!\!\!^\bullet}
\newcommand{\UA}{\Uparrow\!\!\!_\bullet}
\newcommand{\Ua}{\Uparrow}
\newcommand{\Lra}{\Longrightarrow}
\newcommand{\Ra}{\Rightarrow}
\newcommand{\Lla}{\Longleftarrow}
\newcommand{\La}{\Leftarrow}
\newcommand{\nperp}{\perp\!\!\!\!\!\setminus\;\;}
\newcommand{\pq}{\preceq}

\newcommand{\lms}{\longmapsto}
\newcommand{\ms}{\mapsto}
\newcommand{\subseteqnot}{\subseteq\hskip-4 mm_\not\hskip3 mm}

\def\o{{\omega}}

\def\bA{{\bf A}}
\def\bEM{{\bf EM}}
\def\bM{{\bf M}}
\def\bE{{\bf E}}
\def\bN{{\bf N}}
\def\bF{{\bf F}}
\def\bC{{\bf C}}
\def\bI{{\bf I}}
\def\bK{{\bf K}}
\def\bL{{\bf L}}
\def\bT{{\bf T}}
\def\bS{{\bf S}}
\def\bD{{\bf D}}
\def\bB{{\bf B}}
\def\bW{{\bf W}}
\def\bP{{\bf P}}
\def\bX{{\bf X}}
\def\bY{{\bf Y}}
\def\bZ{{\bf Z}}
\def\ba{{\bf a}}
\def\bb{{\bf b}}
\def\bc{{\bf c}}
\def\bd{{\bf d}}
\def\bh{{\bf h}}
\def\bi{{\bf i}}
\def\bj{{\bf j}}
\def\bk{{\bf k}}
\def\bm{{\bf m}}
\def\bn{{\bf n}}
\def\bp{{\bf p}}
\def\bq{{\bf q}}
\def\be{{\bf e}}
\def\br{{\bf r}}
\def\bi{{\bf i}}
\def\bs{{\bf s}}
\def\bt{{\bf t}}
\def\jeden{{\bf 1}}
\def\dwa{{\bf 2}}
\def\trzy{{\bf 3}}

\def\cB{{\cal B}}
\def\cA{{\cal A}}
\def\cC{{\cal C}}
\def\cD{{\cal D}}
\def\cE{{\cal E}}
\def\cEM{{\cal EM}}
\def\cF{{\cal F}}
\def\cG{{\cal G}}
\def\cI{{\cal I}}
\def\cJ{{\cal J}}
\def\cK{{\cal K}}
\def\cL{{\cal L}}
\def\cN{{\cal N}}
\def\cM{{\cal M}}
\def\cO{{\cal O}}
\def\cP{{\cal P}}
\def\cQ{{\cal Q}}
\def\cR{{\cal R}}
\def\cS{{\cal S}}
\def\cT{{\cal T}}
\def\cU{{\cal U}}
\def\cV{{\cal V}}
\def\cW{{\cal W}}
\def\cX{{\cal X}}
\def\cY{{\cal Y}}

%categories

%of functors and monads
%\def\Mnd{{\bf Mnd}}
%\def\MND{{\bf MND}}
%\def\AMnd{{\bf AnMnd}}
%\def\ANMND{{\bf ANMND}}
%\def\An{{\bf An}}
%\def\AN{{\bf AN}}
%\def\Poly{{\bf Poly}}
%\def\SAN{{\bf SAN}}
%\def\San{{\bf San}}
%\def\iSanMnd{_\infty{\bf SanMnd}}
%\def\iSan{_\infty{\bf San}}
%\def\Taut{{\bf Taut}}
%\def\PMnd{{\bf PolyMnd}}
%%\def\RegMnd{{\bf RegMnd}}
%\def\SanMnd{{\bf SanMnd}}
%\def\SANMND{{\bf SANMND}}
%\def\SanMnd{{\bf SanMnd}}
%\def\RiMnd{{\bf RiMnd}}
%\def\End{{\bf End}}
%\def\END{{\bf END}}
%\def\PROP{{\bf PROP}}
%\def\APROP{{\bf AnPROP}}

%of theories
%\def\ET{\bf ET}
%\def\RegET{\bf RegET}
%\def\RET{\bf RegET}
%\def\LrET{\bf LrET}
%\def\RiET{\bf RiET}
%\def\SregET{\bf SregET}
%\def\Cart{\bf Cart}
%\def\wCart{\bf wCart}
%\def\CartMnd{\bf CartMnd}
%\def\wCartMnd{\bf wCartMnd}

%of Lawvere theories
\def\LT{\bf LT}
\def\RegLT{\bf RegLT}
\def\ALT{\bf AnLT}
\def\RiLT{\bf RiLT}

%of Operads
\def\FOp{{\bf FOp}}
\def\RegOp{{\bf RegOp}}
\def\SOp{{\bf SOp}}
\def\Op{{\bf Op}}
\def\RiOp{{\bf RiOp}}
\def\coprodo{{\coprod\!\!\!\!^\ra}}
\def\sqcupo{{\sqcup\!\!^\ra}}

%various other categories and such
%\def\bCat{{{\bf Cat}}}
\def\CAT{{{\bf CAT}}}
\def\Ord{{\bf Ord}}
\def\Card{{\bf Card}}
\def\CAT{{{\bf CAT}}}
\def\MonCat{{{\bf MonCat}}}
\def\Mon{{{\bf Mon}}}
\def\Act{{{\bf Act}}}
\def\mon{{{\bf mon}}}
\def\act{{{\bf act}}}
\def\WMon{{{\bf WMon}}}
\def\WBiMon{{{\bf WBiMon}}}
\def\Cat{{{\bf Cat}}}
\def\Lo{{{\bf Lo}}}
\def\BT{{{\bf BTr}}}

\def\Tl{{{\bf Tl}}}
\def\Ac{{{\bf Ac}}}

\def\F{\mathds{F}}
\def\S{\mathds{S}}
\def\I{\mathds{I}}
\def\B{\mathds{B}}
\def\Lw{\mathds{L}}
\def\T{\mathds{T}}

%functors
\def\V{\mathds{V}}
\def\W{\mathds{W}}
\def\M{\mathds{M}}
\def\A{\mathds{A}}
\def\P{\mathds{P}}
\def\AC{\mathds{AC}}
\def\WAC{\mathds{WAC}}
\def\SM{\mathds{SM}}
\def\WBM{\mathds{WB}i\mathds{M}on}
\def\WCM{\mathds{WC}mon}
\def\WM{\mathds{WM}on}
\def\WSM{\mathds{WSM}on}
\def\N{\mathds{N}}
\def\R{\mathds{R}}
\def\PM{\mathds{PM}}

\def\Vb{\bar{\mathds{V}}}
\def\Wb{\bar{\mathds{W}}}

%\def\Sym{{\cal S}ym}

%leftovers
%\def\P{{\cal P}}
%\def\Q{{\cal Q}}

\pagenumbering{arabic} \setcounter{page}{1}

\title{\bf\Large Weights for Monoids and Actions of Monoids}

\author{{\L}ukasz Sienkiewicz, \\Marek Zawadowski\\ Instytut Matematyki, Uniwersytet Warszawski}
\maketitle
\begin{abstract}
The main objective of the paper is to define the category of monoids as a weighted limit. We also define the category of actions of monoids along the action of a monoidal category as a weighted limit.
\end{abstract}

\section{Introduction}
Weighted limits and colimits provide a uniform way to define many interesting operations on $2$-categories. It has been known for more than 40 years \cite{Law} that the Eilenberg-Moore object for a monad $T$ in a $2$-category $\cK$ is a weighted limit. 

This is an example of a `2-algebraic set' (EM-object) over a 2-dimensional algebraic structure (monad) that can be defined in any sufficiently complete 2-category. One can think that there should be a similar definition of a `2-algebraic set of monoids' over any monoidal category object, another 2-dimensional algebraic structure that can be defined in any sufficiently complete 2-category with finite product. That was a question Bob Pare asked in a personal communication with the second author many years ago, with a further comment `After all, a monoid is a bunch of objects and morphisms satisfying some identities'. The main purpose of this paper is to provide a positive answer to this question by constructing the weight 2-functor $\W$ for the category of monoids. We also construct the weight 2-functor $\mathds{W}_a$ for actions of monoids along an action of monoidal category, to show how one can present many-sorted algebraic structures as weighted limits.

In Sections \ref{sec-monoids} and \ref{sec-actions}  we present the weight 2-functor for categories of monoids ($\W$) and  actions of monoids ($\W_a$), respectively.
Then in section \ref{sec-remarks} we shall comment on these constructions. The paper ends with Appendix containing definitions of monoidal category object, action of a monoidal category object, the category of monoids and the category of actions of monoids.

\subsection*{Notation}

In this paper weighted limits in 2-categories are always meant to be pseudo-limits (i.e., unique up to an iso) and we call them simply (weighted) limits. $\Cat$ is a $2$-category of small\footnote{Note that how `small' are our categories is up to us. The (weighted) limits that we are considering in the paper are always countable.} $2$-categories, functors, and natural transformations. $2\CAT$ is the $3$-category of $2$-categories, i.e., with $2$-categories as $0$-cells, $2$-functors as $1$-cells, $2$-natural transformations as $2$-cells, and $2$-modifications as $3$-cells. Thus $\Cat$ is a $0$-cell of $2\CAT$.
By a $2$-category with finite products we will always mean a $2$-category with finite products of $0$-cells. Let $2\CAT_\times$ be the sub-$3$-category of $2\CAT$ full on $2$-transformations and $2$-modifications, whose $0$-cells are $2$-categories with finite products, and $1$-cells are $2$-functors  preserving
finite products.
$\Mon_{st}(\Cat)$ is a 2-category of monoidal categories, strict monoidal functors, and monoidal natural transformation.
$\Act_{st}(\Cat)$ is a 2-category of actions of monoidal categories, strict morphisms of actions, and transformations of actions.
See Appendix for details. $\Cat$, $\Mon_{st}(\Cat)$, and $\Act_{st}(\Cat)$ are  $0$-cells in $2\CAT_\times$.

$\o$ denotes the set of natural numbers, $n=\{ 0, \ldots, n-1\}$, for $n\in\o$.

\section{The weight 2-functor for the category of monoids}\label{sec-monoids}

Let $M_n$ denote the universe of free monoid on $n$-generators, for $n\in \o+1$. The set $BW_n$ of {\em binary words} on letters $n$  is the free magma on two operations of arity 0 and 2, i.e., it is the least set such that
\begin{enumerate}
  \item $\iota \in BW_n$, where $\iota$ is a distinguished element,
  \item $k\in BW_n$, for $k\in n$,
  \item $(u \diamond v) \in BW_n$, for $u,v\in BW_n$.
\end{enumerate}
We have a type function $ty: BW_\o \ra Lin$ associating to any tree $t$ the linear order  of occurrences of the numbers in $t$, i.e., it is the composition of functions

\[ ty: BW_\o \lra M_\o \lra Lin \]
where the first map is the obvious homomorphism (of magmas) and the second one associates the linear order of occurrences of numbers in the words. For $o\in ty(t)$, we denote by $|o|\in \o$ the number corresponding to the occurrence $o$.

The 2-Lawvere theory $\mathds{M}$ for monoidal categories is a 2-category defined as follows. The objects of $\mathds{M}$ are natural numbers. The morphism in $\mathds{M}$
\[ \vec{u}=(u_0,\ldots, u_{n-1}) : m\ra n\]
is an $n$-tuple of binary word in $BW_m$, i.e., $u_i\in BW_m$, for $i\in n$. The identities are given by
\[ \vec{u} : n\ra n\]
where $u_i=i$, for $i\in n$ and $n\in \o$. The $k$-projections
\[ \pi_k: n\ra 1\]
are given by the binary word $\pi_k=k$, for $k\in n$, , for $n\in \o$.
The composition of morphisms is defined by simultaneous substitution. Given two morphisms
\[  (u_0,\ldots,u_{m-1}) : k\ra m, \hskip 5mm  (v_0,\ldots,v_{n-1}) : m\ra n \]
their composition is given by
\[ (v_0(_0\backslash u_0,\ldots, _{m-1}\backslash u_{m-1}),\ldots, v_{n-1}(_0\backslash u_0,\ldots, _{m-1}\backslash u_{m-1})   ) : k\ra n. \]
and will be shortened to
\[ (v_i(_j\backslash u_j)_{j\in m})_{i\in n} : k\ra n. \]
A unique 2-cell
\[    \vec{u} \Ra \vec{v} : m\ra n \]
exists iff $ty(u_i)=ty(v_i)$, for $i\in n$. Thus hom-sets in $\mathds{M}$ are equivalence relations.
\vskip 2mm

{\em Examples.}
There are morphisms in $\mathds{M}(3,1)$, i.e., 2-cells in $\mathds{M}$,  between $(((1\diamond 2) \diamond (\iota\diamond 0))\diamond 1) $ and     $(((1\diamond 2) \diamond (\iota \diamond(0\diamond 1)))\diamond \iota )$ but not between $(0\diamond (2 \diamond 0))$ and  $(0\diamond (0 \diamond 2))$

The following pictures present binary words as binary trees with labelled leaves

\begin{center} \xext=3500 \yext=550
\begin{picture}(\xext,\yext)(\xoff,\yoff)
%t0:
 \put(-100,0){$t_0=((\iota\diamond 2)\diamond(0\diamond 1))$}
\put(00,400){\line(1,-1){101}}
\put(100,300){\line(1,1){101}}

\put(400,400){\line(1,-1){101}}
\put(500,300){\line(1,1){101}}

\put(100,300){\line(1,-1){200}}
\put(300,100){\line(1,1){200}}

\put(-30,430){$\iota$}
\put(170,430){$2$}
\put(370,430){$0$}
\put(570,430){$1$}

%t1
 \put(950,0){$t_1$}
\put(1000,400){\line(1,-1){300}}
\put(1200,200){\line(1,1){101}}

\put(1200,400){\line(1,-1){101}}
\put(1300,300){\line(1,1){101}}

%\put(100,300){\line(1,-1){200}}
\put(1300,100){\line(1,1){300}}

\put(970,430){$2$}
\put(1170,430){$0$}
\put(1370,430){$\iota$}
\put(1570,430){$1$}

%t2:
 \put(1950,0){$t_2$}
\put(2000,400){\line(1,-1){101}}
\put(2100,300){\line(1,1){101}}

\put(2400,400){\line(1,-1){101}}
\put(2500,300){\line(1,1){101}}

\put(2100,300){\line(1,-1){200}}
\put(2300,100){\line(1,1){200}}

\put(1970,430){$1$}
\put(2170,430){$2$}
\put(2370,430){$0$}
\put(2570,430){$\iota$}

%t3:
 \put(2950,0){$t_3$}
\put(3000,400){\line(1,-1){101}}
\put(3100,300){\line(1,1){101}}

%\put(3400,400){\line(1,-1){100}}
%\put(3500,300){\line(1,1){100}}

\put(3100,300){\line(1,-1){101}}
\put(3200,200){\line(1,1){200}}

\put(2970,430){$1$}
\put(3170,430){$2$}
\put(3370,430){$0$}
%\put(3570,430){$\iota$}

 \end{picture}
\end{center}

\begin{center} \xext=3500 \yext=550
\begin{picture}(\xext,\yext)(\xoff,\yoff)
%t4:
 \put(-50,0){$t_4$}
\put(00,400){\line(1,-1){101}}
\put(100,300){\line(1,1){101}}

\put(400,400){\line(1,-1){101}}
\put(500,300){\line(1,1){101}}

\put(100,300){\line(1,-1){200}}
\put(300,100){\line(1,1){200}}

\put(-30,430){$\iota$}
\put(170,430){$\iota$}
\put(370,430){$1$}
\put(570,430){$1$}

%t5
 \put(950,0){$t_5=((\iota\diamond (1\diamond\iota))\diamond 1)$}
\put(1000,400){\line(1,-1){300}}
\put(1200,200){\line(1,1){101}}

\put(1200,400){\line(1,-1){101}}
\put(1300,300){\line(1,1){101}}

%\put(100,300){\line(1,-1){200}}
\put(1300,100){\line(1,1){300}}

\put(970,430){$\iota$}
\put(1170,430){$1$}
\put(1370,430){$\iota$}
\put(1570,430){$1$}

%t6:
 \put(1950,0){$t_6$}
\put(2000,400){\line(1,-1){101}}
\put(2100,300){\line(1,1){101}}

\put(2400,400){\line(1,-1){101}}
\put(2500,300){\line(1,1){101}}

\put(2100,300){\line(1,-1){200}}
\put(2300,100){\line(1,1){200}}

\put(1970,430){$1$}
\put(2170,430){$2$}
\put(2370,430){$\iota$}
\put(2570,430){$2$}

%t7:
 \put(2950,0){$t_7$}
\put(3000,400){\line(1,-1){101}}
\put(3100,300){\line(1,1){101}}

%\put(3400,400){\line(1,-1){100}}
%\put(3500,300){\line(1,1){100}}

\put(3100,300){\line(1,-1){101}}
\put(3200,200){\line(1,1){200}}

\put(2970,430){$2$}
\put(3170,430){$2$}
\put(3370,430){$1$}
%\put(3570,430){$\iota$}

 \end{picture}
\end{center}
They are 1-cells in $\mathds{M}(3,1)$. There are 2-cells between 1-cells $t_0$ and $t_1$, $t_2$ and $t_3$, $t_4$ and $t_5$, and there are no other 2-cells between these 1-cells.
Moreover, the orders of occurrences can be represented as
\[ ty(t_0)= \lk \lk 0,2 \rk, \lk1,0 \rk, \lk 2,1 \rk \rk, \hskip 10mm ty(t_5)= \lk \lk 0,1 \rk, \lk1,1 \rk \rk, \]
with the first coordinate in the occurrence determining the order and the second the number that is occurring.

\vskip 5mm

The {\em category of words} $CW$ has binary words in $BW_1$ as objects, and a morphism $f:a\ra b$ between two binary words $a,b\in BW_1$ is a monotone function
$f : ty(a) \ra ty(b)$, i.e., between occurrences of $0$'s in $a$ and $b$.

\vskip 5mm

{\em Examples.} Thus $s_0=(0\diamond (\iota \diamond 0))$, $s_1=((((0\diamond (\iota \diamond 0))\diamond \iota)\diamond 0)\diamond 0)$ are objects of $CW$.
The morphism from $s_0$ to $s_1$ in $CW$ is monotone function
\[ f :  \lk \lk 0,0 \rk, \lk1,0 \rk \rk \lra  \lk \lk 0,0 \rk, \lk1,0 \rk,\lk 2,0 \rk, \lk3,0 \rk \rk \]
between linear orders.

\vskip 5mm

The weight 2-functor
\[ \mathds{W}: \mathds{M} \lra \Cat \]
is defined as follows. For $m\in \o$
\[ \mathds{W}(m) = CW^m. \]
For $\vec{u} : m\ra n \in \mathds{M}$, the functor
\[ \mathds{W}(\vec{u}) : CW^m \ra CW^n,  \]
is given on object $\vec{a}=(a_0, \ldots, a_{m-1})\in CW^m$ by
\[ \mathds{W}(\vec{u})(\vec{a}) = (u_i(j\backslash a_j)_{j\in m})_{i\in n}, \]
i.e., we substitute $m$ objects from $BW_1$ into $n$ objects from $BW_m$ and we obtain $n$ objects in $BW_1$.

On morphism $\vec{f}: \vec{a}\ra \vec{b}\in CW^m$ the functor $\mathds{W}$ is given by
\[ \mathds{W}(\vec{u})(\vec{f}) =  (\coprodo_{o\in ty(u_i)} f_{|o|})_{i\in n}, \]
where $\coprodo_{x\in X} Y_x$ is the ordered sum of posets $Y_x$ indexed by a poset $X$. In particular, if $X$ and $Y_x$ are linear orders, the sum is a linear order as well. $X\sqcupo Y$ is an ordered sum of two posets.

To a 2-cell
\[ \alpha :  \vec{u} \Ra \vec{v} : m\ra n \]
the 2-functor $\mathds{W}$ associates a natural transformation
\[ \mathds{W}(\alpha) : \mathds{W}(\vec{u}) \ra   \mathds{W}(\vec v), \]
so that its $i$-th component at $\vec{a}\in CW^m$ is

\[  {\mathds{W}(\alpha)_{\vec{a},i}} = \coprodo_{o\in ty(u_i) } id_{ty(a_{|o|})} : (u_i(j\backslash a_j)_{j\in m}) \ra (v_i(j\backslash a_j)_{j\in m}, \]
for $i\in n$.

The fact that  $\mathds{W}(\alpha)$ is indeed a natural transformation follows from the commutativity of the following diagram:
\begin{center} \xext=2200 \yext=1000
\begin{picture}(\xext,\yext)(\xoff,\yoff)
\setsqparms[1`1`1`1;2200`800]
\putsquare(0,100)[(u_i(j\backslash a_j)_{j\in m})_{i\in n}`(v_i(j\backslash a_j)_{j\in m})_{i\in n}`(u_i(j\backslash b_j)_{j\in m})_{i\in n}`(v_i(j\backslash b_j)_{j\in m})_{i\in n};{\mathds{W}(\alpha)_{\vec{a}}}=(\coprodo_{o\in ty(u_i) } id_{ty(a_{|o|})})_{i\in n}`{\mathds{W}(\vec{u})(\vec{f}})=(\coprodo_{o\in ty(v_i) } f_{|o|})_{i\in n}`(\coprodo_{o\in ty(v_i) } f_{|o|})_{i\in n}={\mathds{W}(\vec{v})(\vec{f}})`{\mathds{W}(\alpha)_{\vec{b}}}=(\coprodo_{o\in ty(u_i) } id_{ty(b_{|o|})})_{i\in n}]
\end{picture}
\end{center}
in $\mathds{W}(n)$, which is an ordered sum of the squares

\begin{center} \xext=800 \yext=650
\begin{picture}(\xext,\yext)(\xoff,\yoff)
\setsqparms[1`1`1`1;800`500]
\putsquare(0,50)[a_{|o|}`a_{|o|}`b_{|o|}`b_{|o|};id_{ty(a_{|o|})}`f_{|o|}` f_{|o|}`id_{ty(b_{|o|})}]
\end{picture}
\end{center}
for $o\in ty(u_i)=ty(v_i)$ and $i\in n$.

\vskip 2mm

\begin{theorem}
We have a equivalence of 2-categories \[ \zeta:\Mon_{st}(\Cat)\lra 2\CAT_\times(\M,\Cat)\] such that the triangle
\begin{center} \xext=600 \yext=650
\begin{picture}(\xext,\yext)(\xoff,\yoff)
 \settriparms[1`1`1;600]
  \putVtriangle(0,0)[\Mon_{st}(\Cat)`2\CAT_\times(\M,\Cat)`\Cat.;\zeta`\mon`Lim_\W]
\end{picture}
\end{center}
commutes up to a 2-equivalence.
In other words, $\M$ is the 2-Lawvere theory for monoidal categories and $\mathds{W} : \mathds{M} \lra \Cat$ is the weight 2-functor for the category of monoids.
\end{theorem}

{\em Proof.}
$\mathds{M}$ has finite products given by the sum: $m\times n = m+n$.

To see the first claim, we shall define the 2-functor
\begin{center}
\xext=600 \yext=200
\begin{picture}(\xext,\yext)(\xoff,\yoff)
 \putmorphism(0,50)(1,0)[\Mon_{st}(\Cat)`2\CAT_\times(\M,\Cat)`\zeta]{1200}{1}a
\end{picture}
\end{center}
from the 2-category of monoidal categories to the 2-category for finite product preserving 2-functors $F:\cM \ra \Cat$, natural transformations and modifications between them.

Let $(\cM,\otimes,I,\alpha,\lambda,\rho)$ be a monoidal category. We define a 2-functor

\[ \overline{\cM}=\zeta(\cM,\otimes,I,\alpha,\lambda,\rho) : \mathds{M} \lra \Cat \]
as follows. For $n\in \o$
\[ \overline{\cM}(n) = \cM^n. \]
for $s: n \ra 1$ in $\mathds{M}$

 \[  \overline{\cM}(s) \;= \; \left\{ \begin{array}{ll}
		I     & \mbox{ if } s=\iota, \\
        \pi_k & \mbox{ if } s=k\in n, \\
        \otimes \circ \lk \overline{\cM}(s_1),\overline{\cM}(s_2) \rk &  \mbox{ if } s=(s_1\diamond s_2). \\
                                    \end{array}
			    \right. \]
where $I : \cM^n\ra \cM$ is the constant functor with value $I$, and $\pi_k: \cM^n\ra \cM$ is the $k$-th projection. As $\overline{\cM}$ is supposed to preserve finite products, it is enough to define it on morphisms with codomain $1$.

If $\sigma : s \Ra t : n \ra 1$ is a 2-cell in $\mathds{M}$, then $\overline{\cM}(\sigma) : \overline{\cM}(s) \Ra \overline{\cM}(t):\overline{\cM}(n) \Ra \overline{\cM}(1) $ is the unique formal (i.e., built from $\alpha$, $\lambda$, $\rho$) natural transformation between functors $\overline{s}$ and $\overline{t}$. Such a natural transformation exists and is unique by MacLane's coherence theorem for monoidal categories.

In order to show that $\mathds{W}: \mathds{M} \ra \Cat$ is indeed the weight 2-functor for the category of monoids, we shall construct the universal $\mathds{W}$-weighted cone $\tau$ over $\overline{\cM}$ with the vertex being the category of monoids $\mon(\cM,\ldots)=\mon(\cM,\otimes,I,\alpha,\lambda,\rho)$. As $\mathds{W}$ preserves finite products, it is enough to define the projections  $\tau_{(-)}: \mon(\cM)\lra \cM$ indexed by the objects and morphisms of category $\mathds{W}(1)$.

For object $a\in  \mathds{W}(1)$ and monoid $(M,m,i)$, we put
 \[  \tau_{a}(M,m,i) \;= \; \left\{ \begin{array}{ll}
		I     & \mbox{ if } a=\iota, \\
        M & \mbox{ if } a=0, \\
        \tau_{a_1}(M,m,i)\otimes \tau_{a_2}(M,m,i)  &  \mbox{ if } a=(a_1\diamond a_2), \\
                                    \end{array}
			    \right. \]
and for homomorphism $h:(M,m,i)\ra (M',m',i')$
 \[  \tau_{a}(h) \;= \; \left\{ \begin{array}{ll}
		1_I     & \mbox{ if } a=\iota, \\
        h & \mbox{ if } a=0, \\
        \tau_{a_1}(h)\otimes \tau_{a_2}(h) &  \mbox{ if } a=(a_1\diamond a_2). \\
                                    \end{array}
			    \right. \]

We still need to define $\tau$ on morphisms of $\mathds{W}(1)$. We shall do it for the (unique) morphisms of the form $f_a: a\ra 0$. The unique extension to all the morphisms in $\mathds{W}(1)$ is again due to the  MacLane's coherence theorem. We define

  \[  \tau_{f_a}(M,m,i) \;= \; \left\{ \begin{array}{ll}
		i     & \mbox{ if } a=\iota, \\
        1_M & \mbox{ if } s=0, \\
        m \circ  ( \tau_{f_{a_1}}(M,m,i)\otimes\tau_{f_{a_2}}(M,m,i)) &  \mbox{ if } a=(a_1\diamond a_2). \\
                                    \end{array}
			    \right. \]
The remaining details of the construction and verifications that $\tau$ is indeed a universal $\W$-weighted cone are left for the reader.
$\boxempty$

\section{The weight 2-functor for actions of monoids}\label{sec-actions}

In this section we construct the weight 2-functor for actions of monoids along an action of a monoidal category. As the Lawvere theory of actions of monoidal categories has two types,  we shall be using two (disjoint) sets of natural numbers: the usual one $\o$, and the set of `underlined' natural numbers $\underline{\o}=\{ \underline{n}: n\in \o\}$.

The set $BW_{*,n,n'}$ of {\em pointed binary words} on letters $n\in\o+1$ and $n'\in \o+1$  is the least set such that
\begin{enumerate}
  \item $\underline{k}\in BW_{*,n,n'}$, for $k\in n'$,
  \item $(u \diamond s) \in BW_{*,n,n'}$, for $u\in BW_n$ and $s\in BW_{*,n,n'}$.
\end{enumerate}

We have a type function $ty_*: BW_{*,\o,\o} \ra Lin_*$ associating to any pointed binary word $t$ the linear order with the right end-point of occurrences of the numbers in $t$. Note that in pointed binary words there is always one occurrence of an underlined number and it is the last occurrence in the word; so it is the top end-point in the linear order of occurrences. If $o\in ty_*(t)$, then $|o|$ denotes the number or underlined number corresponding to the occurrence $o$.

The 2-Lawvere (finite product) theory $\mathds{A}$, for actions of monoidal category,  is a 2-category defined as follows. The objects of $\mathds{A}$ are pairs of natural numbers. The 1-cells in $\mathds{A}$
\[ (\vec{u},\vec{s})=(u_0,\ldots, u_{n-1};s_0,\ldots, s_{n'-1}) : (m,m')\lra (n,n')\]
is an $n$-tuple of binary words in $BW_m$, i.e., $u_i\in BW_m$, for $i\in n$ and $n'$-tuple of pointed binary words in $BW_{*,m,m'}$, i.e., $s_j\in BW_{*,m,m'}$, for $i\in n'$.
When convenient, we shorten the notation for morphisms to $(u_i;s_{i'})_{i\in n, i'\in n'}$ or even to $(\vec{u};\vec{s})$.

The identity is given by
\[ (u_0,\ldots, u_{n-1};s_0,\ldots, s_{n'-1} ) : (n,n')\lra (n,n',)\]
where $u_i=i$, for $i\in n$, and $s_{i'}=\underline{i'}$, for $i'\in n'$. The $k$-th projection is
\[ \pi_k=(u; ) : (n,n')\lra (1,0)\]
where $u=k$, for $k\in n$, and the $\underline{k}$-th projection is
\[ \pi_{\underline{k}}=(;s) : (n,n')\lra (0,1)\]
where $s=\underline{k}$, for $k\in n'$.

The composition of morphisms is defined by simultaneous substitution. Given two morphisms
\[  (u_0,\ldots,u_{m-1};s_0,\ldots, s_{m'-1} ) : (k,k')\lra (m,m'), \hskip 5mm   (v_0,\ldots,v_{n-1};t_0,\ldots, t_{n'-1}) : (m,m')\lra (n,n') \]
their composition is given by
\[ (v_i(_j\backslash u_j)_{_j\in m};t_{i'}(_j\backslash u_j  ,_{\underline{j'}}\backslash s_{j'})_{j\in m, j'\in m'} )_{i\in n, i'\in n'} : (k,k')\lra (n,n'). \]
A unique 2-cell
\[   (u_i;s_j)_{i\in n, j\in n'} \Ra (v_i;t_j)_{i\in n, j\in n'} : (m,m')\lra (n,m') \]
exists iff $ty(u_i)=ty(v_i)$, for $i\in n$, and $ty_*(s_j)=ty_*(t_j)$, for $j\in n'$.

\vskip 2mm

The {\em category of pointed words} $CW_*$ has $BW_{*,1,1}$ as the set of objects, and a morphism $f:a\ra b$ between two pointed binary words $a,b\in BW_{*,1,1}$ is monotone top end-point preserving function $f : ty(a) \ra ty(b)$, i.e., between occurrences of $0$'s and $\underline{0}$ (occurrence of $0$ can be sent to the occurrence $\underline{0}$ but not vice versa).

The weight 2-functor
\[ \mathds{W}_a: \mathds{A} \lra \Cat \]
is defined as follows.
For $m,m'\in \o$
\[ \mathds{W}_a(m,m') = CW^m \times CW_*^{m'}. \]
For $(\vec{u};\vec{s}) : (m,m')\lra (n,n') \in \mathds{A}$, the functor
\[ \mathds{W}_a(\vec{u};\vec{s}) : \mathds{W}_a(m,m')\lra \mathds{W}_a(n,n'),  \]
is given, for $(\vec{a};\vec{b})\in \mathds{W}_a(m,m')$, by
\[ \mathds{W}_a(\vec{u};\vec{s})(\vec{a};\vec{b}) = (u_i(_j\backslash a_j)_{_j\in m};s_{i'}(_j\backslash a_j;_{\underline{j'}}\backslash b_{j'})_{j\in m, j'\in m'}  )_{i\in n, i'\in n'} \]
and, on morphism $(\vec{f},\vec{g}): (\vec{a};\vec{b})\ra (\vec{a}';\vec{b}')\in \mathds{W}_a(m,m')$, is given by
\[ \mathds{W}_a(\vec{u};\vec{s})(\vec{f};\vec{g}) =  (\coprodo_{o\in ty(u_i)} f_{|o|}; \coprodo_{o\in ty(s_{i'}),\; |o|\in \o} f_{|o|} \sqcupo  \coprodo_{o\in ty(s_{i'}),\; |o|\in \underline{\o}} g_{|o|}   )_{i\in n, i'\in n'}. \]
Note that $\mathds{W}_a(\vec{u};\vec{s})(\vec{f};\vec{g})$ is well defined as the last summand of the right sum contains exactly one element. To a 2-cell
\[ \alpha :  (\vec{u};\vec{s}) \Ra (\vec{v};\vec{t}) : (m,m')\lra (n,n') \]
the 2-functor $\mathds{W}_a$ associates a natural transformation
\[ \mathds{W}_a(\alpha) : \mathds{W}_a(\vec{u};\vec{s}) \lra   \mathds{W}_a(\vec{v};\vec{t}): \mathds{W}_a(m,m') \lra   \mathds{W}_a(n,n'), \]
so that its $i$-th component at $(\vec{a},\vec{b})\in \mathds{W}_a(m,m')$ is

\[  {\mathds{W}_a(\alpha)_{\vec{a},\vec{b},i}} = \coprodo_{o\in ty(u_i) } id_{ty(a_{|o|})} : (u_i(j\backslash a_j)_{j\in m}) \lra (v_i(j\backslash a_j)_{j\in m}, \]
for $i\in n$, and its $\underline{i}'$-th component is
\[  {\mathds{W}_a(\alpha)_{\vec{a},\vec{b},i'}} = \coprodo_{o\in ty(s_{i'}),\; |o|=0 } id_{ty(a_{|o|})} \sqcupo \coprodo_{o\in ty(s_{i'}),\; |o|=\underline{0} } id_{ty(b_{|o|})} : (s_{i'}(_j\backslash a_j;_{j'}\backslash b_{j'})_{j\in m}) \lra (t_{i'}(j\backslash a_j;_{j'}\backslash b_{j'})_{j\in m, j'\in m'}, \]
for $i'\in n'$.

The fact that  $\mathds{W}_a(\alpha)$ is indeed a natural transformation follows from the commutativity of the following diagram in $\Cat$:

\begin{center} \xext=2200 \yext=1000
\begin{picture}(\xext,\yext)(\xoff,\yoff)
\setsqparms[1`1`1`1;2200`800]
\putsquare(50,100)[\mathds{W}_a({\vec{u};\vec{s}})({\vec{a};\vec{b}})`\mathds{W}_a{(\vec{v};\vec{t})}({\vec{a};\vec{b}})`\mathds{W}_a{(\vec{u};\vec{s})}({\vec{a}';\vec{b}'})`
\mathds{W}_a{(\vec{v};\vec{t})}({\vec{a}';\vec{b}'});\mathds{W}_a(\alpha)_{({\vec{a};\vec{b}})}`\mathds{W}_a({\vec{u};\vec{s}})(\vec{f};\vec{g})`
\mathds{W}_a{(\vec{v};\vec{t})}(\vec{f};\vec{g})`\mathds{W}_a(\alpha)_{({\vec{a}';\vec{b}'})}]
\end{picture}
\end{center}
for any morphism $(\vec{f};\vec{g}): (\vec{a};\vec{b})\lra (\vec{a}';\vec{b}')$ in $\mathds{W}_a(n,n')$. And this can be checked in a similar way as in the case of the weight 2-functor for monoids $\mathds{W}$.

\begin{theorem}
We have a equivalence of 2-categories \[ \zeta_a:\Act_{st}(\Cat)\lra 2\CAT_\times(\mathds{A},\Cat)\] such that the triangle
\begin{center} \xext=600 \yext=650
\begin{picture}(\xext,\yext)(\xoff,\yoff)
 \settriparms[1`1`1;600]
  \putVtriangle(0,0)[\Act_{st}(\Cat)`2\CAT_\times(\mathds{A},\Cat)`\Cat.;\zeta_a`\act`Lim_{\W_a}]
\end{picture}
\end{center}
commutes up to a 2-isomorphism.
In other words, $\mathds{A}$ is the (2-sorted) 2-Lawvere theory for actions of monoidal categories and $\mathds{\W_a} : \mathds{A} \lra \Cat$ is the weight 2-functor for the category of actions of monoids along an action of a monoidal category.
\end{theorem}

{\em Proof.} $\A$ has finite products given by the sum: $(n,m)\times (n',m') = (n+n',m+m')$.

Let $\cA=(\cM,\otimes,I,\alpha,\lambda,\rho,\cX,\star, \psi,\bar{\psi})$ be an action a monoidal category $(\cM,\otimes,I,\alpha,\lambda,\rho)$  on a category $\cX$. As before, we shall describe the universal $\W_a$-cone with the vertex being the category of actions  $\act(\cA)=\act(\cM,\otimes,I,\alpha,\lambda,\rho,\cX,\star, \psi,\bar{\psi})$.

We define a 2-functor

\[ \overline{\cA}=\zeta_a(\cA) : \mathds{A} \lra \Cat \]
as follows. For $n,n'\in \o$
\[ \overline{\cA}(n,n') = \cM^n\times \cX^{n'}. \]
For a morphism $(u,\emptyset): (n,n') \ra (1,0)$ in $\mathds{A}$, we put
\[ \overline{\cA}(u,\emptyset)=\overline{\cM}(u)\circ \pi_{\cM} : \cM^n\times \cX^{n'}\lra \cM \]
where $\pi_{\cM}: \cM^n\times \cX^{n'}\ra  \cM^n$ is the obvious projection, and for a morphism $(\emptyset,s): (n,n') \ra (0,1)$ in $\mathds{A}$, we put
\[ \overline{\cA}(\emptyset,s)= : \cM^n\times \cX^{n'}\lra  \cX \]
to be
 \[  \overline{\cA}(\emptyset,s) \;= \; \left\{ \begin{array}{ll}
		\pi_{\underline{k}}     & \mbox{ if } s=\underline{k}, \\
            &\\
         (\overline{\cM}(u')\circ \pi_{\cM})\star \overline{A}(\emptyset,s') &  \mbox{ if } s=(u'\diamond s'), \\
                                    \end{array}
			    \right. \]

As $\overline{\cA}$ is supposed to preserve finite products, it is enough to define it on morphisms with codomain $(1,0)$ and $(0,1)$.

If $\sigma : (u,s) \Ra (v,t) : (n ,n')\ra (1,1)$ is a 2-cell in $\mathds{A}$, then $\overline{\cA}(\sigma) : \overline{\cA}(u,s) \Ra \overline{\cA}(v,t)$ is the unique formal (i.e., built from $\alpha$, $\lambda$, $\rho$, $\psi$, $\bar{\psi}$) natural transformation between these functors. Such a natural transformation exists and is unique, again, by MacLane's coherence theorem for monoidal categories.

In order to show that $\mathds{W}_a: \mathds{A} \ra \Cat$ is indeed the weight 2-functor for the category of monoids, we shall construct the universal $\mathds{W}_a$-weighted cone $\kappa$ over $\overline{\cA}$ with the vertex being the category of actions of monoids along action of monoidal category  $\act(\cA)=\act(\cM,\otimes,I,\alpha,\lambda,\rho,\cX,\star, \psi,\bar{\psi})$. As $\mathds{W}_a$ preserves finite products, it is enough to define the projections
\[ \kappa_{(-)}: \act(\cA)\lra \cM\hskip 5mm {\rm and}\hskip 5mm \kappa_{(-)}: \act(\cA)\lra \cX \]
indexed by the objects and morphisms of categories $\mathds{W}_a(1,0)$ and $\mathds{W}_a(0,1)$, respectively. Note that the projections $\kappa_{a,\emptyset}$ and $\kappa_{f_a,\emptyset}$ composed with $\pi_\cM$ are the projections $\tau_{a}$ $\tau_{f_a}$, respectively, from the universal $\mathds{W}$-cone for monoids $\tau$, defined in the previous proof.

For object $(a,\emptyset)\in  \mathds{W}_a(1,0)$ and action $(M,m,i,X,\alpha)$, we put
\[ \kappa_{a,\emptyset}(M,m,i,X,\alpha) =\tau_a((M,m,i), \]
and for and for morphism of actions $(h,k):(M,X,m,i,\alpha)\ra (M',X',m',i',\alpha')$
\[ \kappa_{a,\emptyset}(h,k) =\tau_a(h). \]
To define $\kappa$ on morphisms of $\mathds{W}_{a}(1,0)$ we again use projection $\tau$. We shall define the naural transformation $\kappa_{(f_a;\emptyset}$  for the (unique) morphisms of the form $(f_a,\emptyset): (a,\emptyset)\ra (0,\emptyset)$ in $\mathds{W}_{a}(1,0)$. The unique extension to all the morphisms in $\mathds{W}_a(1,\emptyset)$ is again due to the  MacLane's coherence theorem. We put
\[ \kappa_{f_a,\emptyset}(M,X,m,i,\alpha) = \tau_{f_a}(M,m,i). \]

For object $(\emptyset,b)\in  \mathds{W}_a(0,1)$ and action $(M,m,i,X,\alpha)$, we put

 \[  \kappa_{\emptyset,b}(M,m,i,X,\alpha) \;= \; \left\{ \begin{array}{ll}
        X & \mbox{ if } b=\underline{0}, \\
        \tau_a(M,m,i)\star \kappa_{\emptyset,b'}(M,m,i,X,\alpha) &  \mbox{ if } b=(a\diamond b'), \\
                                    \end{array}
			    \right. \]
and for morphism of actions $(h,k):(M,X,m,i,\alpha)\ra (M',X',m',i',\alpha')$

 \[  \kappa_{\emptyset,b}(h,k) \;= \; \left\{ \begin{array}{ll}
        k & \mbox{ if } b=\underline{0}, \\
        \tau_{a}((h)\star\kappa_{\emptyset,b'}(h,k)  &  \mbox{ if } b=(a\diamond b'). \\
                                    \end{array}
			    \right. \]

We still need to define $\kappa$ on morphisms of $\mathds{W}_{\emptyset,b}(\emptyset,1)$. We shall do it for the (unique) morphisms of the form $(\emptyset,g_b): (\emptyset,b)\ra (\emptyset,\underline{0})$. The unique extension to all the morphisms in $\mathds{W}_a(\emptyset,1)$ is again due to the  MacLane's coherence theorem. We define

  \[ \kappa_{\emptyset,g_b}(M,m,i,i,\alpha) \;= \; \left\{ \begin{array}{ll}
        1_X & \mbox{ if } b=\underline{0}, \\
        \alpha \circ  (\tau_{f_{a_1}}(M,m,i)\star\kappa_{1,f_{a},g_{b'}}(M,m,i,\alpha) ) &  \mbox{ if } b=(a\diamond b'). \\
                                    \end{array}
			    \right. \]
The remaining details of the construction and verifications that $\kappa$ is indeed a universal $\W_a$-weighted cone are left for the reader.
$\boxempty$

\section{Concluding remarks}\label{sec-remarks}

\subsection*{Extension of the context.}
To make is simple the paper is presented in the essentially simplest context in which considerations of these weights is meaningful. Now we shall comment on possible extension of the context.

We know that in $\Cat$ not only strict but also lax monoidal functors between monoidal categories induce functors between categories of monoids. This is due to the fact even if the composition $H(\tau)$ of a (strict) universal cone $\tau : \mon(\cM,\otimes, I, \alpha, \lambda,\rho) \dot{\lra} \cM$ composed with a lax monoidal functor
\[ (F,\overline{\varphi},\varphi): (\cM,\otimes, I, \alpha, \lambda,\rho)\lra (\cM',\otimes', I', \alpha', \lambda',\rho')\]
is not a strict cone (i.e., the commutations of 1-cells holds up to a non-invertible 2-cells) we can `strictify' such a cone to a cone $\sigma: \mon(\cM,\otimes, I, \alpha, \lambda,\rho \dot{\lra} \cM'$ that is strict and hence it induces a functor between categories of monoids by the universal property of the category of monoids $\mon(\cM',\otimes', I', \alpha', \lambda',\rho')$. The strictified cone $\sigma$ can be defined as follows. For object $a\in  \mathds{W}(1)$ and monoid $(M,m,i)$, we put
 \[  \sigma_{a}(M,m,i) \;= \; \left\{ \begin{array}{ll}
		I'     & \mbox{ if } a=\iota, \\
        F(M) & \mbox{ if } a=0, \\
        \sigma_{a_1}(M,m,i)\otimes \sigma_{a_2}(M,m,i)  &  \mbox{ if } a=(a_1\diamond a_2), \\
                                    \end{array}
			    \right. \]
and for homomorphism $h:(M,m,i)\ra (M',m',i')$
 \[  \sigma_{a}(h) \;= \; \left\{ \begin{array}{ll}
		1_{I'}     & \mbox{ if } a=\iota, \\
        F(h) & \mbox{ if } a=0, \\
        \sigma_{a_1}(h)\otimes \sigma_{a_2}(h) &  \mbox{ if } a=(a_1\diamond a_2). \\
                                    \end{array}
			    \right. \]
For the (unique) morphisms of the form $f_a: a\ra 0$, we define

  \[  \sigma_{f_a}(M,m,i) \;= \; \left\{ \begin{array}{ll}
		F(i)\circ \overline{\varphi}:I'\ra H(M)     & \mbox{ if } a=\iota, \\
        1_M & \mbox{ if } s=0, \\
        H(m) \circ  ( \sigma_{f_{a_1}}(M,m,i)\otimes\sigma_{f_{a_2}}(M,m,i))\circ \varphi_{\sigma_{a_1}(M,m,i), \sigma_{a_2}(M,m,i)} &  \mbox{ if } a=(a_1\diamond a_2). \\
                                    \end{array}
			    \right. \]
The remaining details of the construction and verifications that $\sigma$ is indeed a (strict) universal $\W$-weighted cone over $(\cM,\otimes, I, \alpha, \lambda,\rho)$ is a routine check.

Once we have the weight 2-functor for categories of monoids, we can use it do define the category of monoids objects over any monoidal object in any 2-category with finite products. As the strictification described above can be defined using the universal properties of finite products, it is still true that a lax monoidal morphism of monoidal category objects in any 2-category with finite products induces a 1-cell between objects of monoids.

\subsection*{Algebra needs coalgebra}

Affine algebraic sets over (set-based) algebraic structures (like rings, fields, groups, module etc.) can be defined as limits on the diagrams that involve finite power of the universe and some definable (polynomial) functions between them. Typically, these limits come from finite sets of equations
\[ f(x_1, \ldots, x_n) = g(x_1, \ldots, x_n) \]
but we do not make any restrictions on the variables that occur on both sides of the equations so that we can consider equations like
\[ f(x,x,y)=g(x,y,y,y) \]
that use the same variable more than once, not necessarily the same number of times on each side (thus using diagonals), and we can also have equations
\[ m(x,x) = e \]
that have different variables occurring on different sides of the equation (thus using projections).

The limits giving rise to such algebraic sets can be chosen canonically, if we allow weights in their definitions. We shall `prove it' by an example. Let $A$ be a commutative ring in a complete category $\cE$. Then the equation $x^2=y^3$ defines a subobject $Z$ of the square of the universe of $A$ (also denoted by $A$). In the internal language it can be expressed as
\[ Z= \{ \lk a,b \rk \in A^2 | a^2=b^3   \} \]
Let $B=\cF[x,y]_{/x^2-y^3}$ be the free commutative ring (in $Set$) on two generators $x$ and $y$ divided by the equation $x^2=y^3$, and $\Lw_{cring}$ be the Lawvere theory for commutative rings. We have finite products preserving functors
\[ \bar{A}:\Lw_{cring}\ra \cE ,\hskip 1cm \bar{B} : \Lw_{cring}\ra Set \]
corresponding to the rings $A$ and $B$. Then, the set $Z$ is the weighted limit $Lim_{\bar{B}}\bar{A}$. We show that it is the case if $\cE$ is the category of sets  $Set$. We have a sequence of isomorphisms
\[ Z = \{ \lk a,b \rk \in A^2 | a^2=b^3   \} \cong \]
\[ Hom(B,A) \cong \]
\[ Nat(\bar{B},\bar{A}) \cong \]
\[  Nat(\bar{B},Set(1,\bar{A}(-))) \cong \]
\[ Set(1,Lim_{\bar{B}}\bar{A})\cong \]
\[ Lim_{\bar{B}}\bar{A} \]
where $Hom$ is the hom-set in the category of commutative rings.

Note that we also have a PROP\footnote{PROP is a strict symmetric monoidal category whose objects are natural numbers such that $I=0$ and $n\otimes m=n+m$.} for commutative rings $\P_{cring}$ and hence symmetric monoidal functors
\[ \widetilde{A}:\P_{cring}\ra \cE ,\hskip 1cm  \widetilde{B} : \P_{cring}\ra Set \]
corresponding to rings $A$ and $B$. The monoidal structure considered on both $\cE$ and $Set$ is the finite product structure.  However, it is not the case that $Lim_{\widetilde{B}}\widetilde{A}$ is isomorphic to the object $Z$ (even if $\cE$ is $Set$), as natural transformations from $\widetilde{B}$ to $\widetilde{A}$ do not correspond to homomorphisms from $B$ to $A$ in this case. The reason for this is that $\P_{cring}$ does not have projections and diagonals, a piece of coalgebra which was vital in the former argument.

In other words, to define the usual algebraic sets we use the coalgebra structure on this set with respect to the tensor being the usual cartesian product. This comonoid structure is usually not mentioned for good reasons: it is unique, if our tensor is the binary product. However, if we replace the product by some other tensor, we need to specify the comonoid structure separately, if we want to use it. In this sense to do algebra we need to use a bit of coalgebra.  This must be taken into account when we define 2-algebraic structures, as well.

\subsection*{2-algebra needs 2-coalgebra}

As we already learned from the previous discussion, it is not necessarily true that if we can identify internally an algebraic concept (ring), then we will be able to derive all the `algebraic sets' related to it (set of solutions of equations). In fact, now we need to talk about `2-algebraic sets' as the derived concepts will be categories or even $0$-cells in a $2$-category. The category $s\Delta$ is not even a $2-PROP$, (i.e., a strict symmetric monoidal 2-category whose $0$-cells are natural numbers such that $I=0$ and $n\otimes m=n+m$) but the Eilenberg-Moore object has sufficiently simple structure that we are able to get it as a weighted limit from the functor with domain $s\Delta$. Thus in this case no coalgebra is needed.

If we want to internalize the notion of a monoidal category, we have to have, in our ambient $2$-category $\cK$, finite products of $0$-cells or at least a 2-monoidal structure. Then we can easily define a $2$-category $\PM$ which is a $2-PROP$ for monoidal categories, i.e., with the property that if $\cK$ is a $2$-category with finite products, then the $2$-monoidal functors from $\PM$ to $\cK$ correspond to  monoidal categories in $\cK$. However, if we want to derive a `2-algebraic set' of monoids from a monoidal category, the weighted limit of a monoidal $2$-functor from $\PM$ is not enough. This is because when we look at the structure maps of monoids
\[ m :M\otimes M \ra m\hskip 1cm i:I \ra m\]
as `some kind of equalities', they are not linear-regular (cf. \cite{SZ1}, \cite{SZ2}) as in the left `equation' a variable is repeated and on the right a variable is dropped. Thus this uses full force of the equational logic, not just the linear-regular part. Therefore to define internally the object of monoids either we need the internal version of the notion of a bi-monoidal category (by this we mean the categorification of the notion of a bi-monoid) or we need to define the internal notion of a monoidal category on the basis of finite product, i.e., not using $2-PROP$'s but Lawvere 2-theories. In this paper, we had follow the latter approach but a further extension is still possible using the former.

\section{Appendix.}

\subsection{Monoidal category and category of monoids}\label{subsec-monoidal-obj}

A {\em monoidal category} $(\cM,\otimes,I,\alpha,\lambda,\rho)$ consists of
\begin{enumerate}
  \item a category $\cM$,
  \item a functor $\otimes :\cM\times \cM\ra \cM$, and object $I \in \cM$,
  \item three natural isomorphisms
  \[ \alpha :\otimes \circ (1_{\cM}\times \otimes)\lra \otimes \circ(\otimes\times 1_\cM) : \cM\times \cM\times \cM\ra \cM,\]
   \[ \lambda : \otimes\circ \lk I, 1_\cM\rk\ra 1_\cM : \cM\ra \cM, \hskip 5mm \rho: \otimes\circ \lk 1_\cM,I \rk \ra 1_\cM : \cM\ra \cM\]
\end{enumerate}
such that, for objects $M_1$, $M_2$, $M_3$, $M_4$  in $\cM$, the following two diagrams in $\cM$
\begin{center} \xext=2300 \yext=1100
\begin{picture}(\xext,\yext)(\xoff,\yoff)
\put(500,1020){$M_1\otimes({M_2}\otimes ({M_3}\otimes {M_4}))$}
      \put(900,970){\vector(-3,-2){360}}
      \put(50,870){$_{1_{M_1}\otimes \alpha_{{M_2},{M_3},{M_4}}}$}

      \put(1300,970){\vector(3,-2){360}}
      \put(1520,870){$_{\alpha_{M_1,{M_2},{M_3}\otimes {M_4}}}$}

\put(-200,620){$M_1\otimes(({M_2}\otimes {M_3})\otimes {M_4})$}
\put(1450,620){$(M_1\otimes {M_2})\otimes ({M_3}\otimes {M_4})$}
      \put(500,570){\vector(1,-2){160}}
      \put(30,360){$_{\alpha_{M_1,{M_2}\otimes {M_3},{M_4}}}$}

      \put(1700,570){\vector(-1,-2){160}}
      \put(1640,360){$_{\alpha_{M_1\otimes {M_2},{M_3},{M_4}}}$}

\put(-80,120){$(M_1\otimes({M_2}\otimes {M_3}))\otimes {M_4}$}
\put(1400,120){$((M_1\otimes {M_2})\otimes {M_3})\otimes {M_4}$}

  \put(1050,120){\vector(1,0){300}}
      \put(950,0){$_{\alpha_{M_1,{M_2},{M_3}}\otimes 1_{M_4}}$}
       \put(-800,400){$\bf MC1$}
\end{picture}
\end{center}
and
  \begin{center}\xext=1000 \yext=700
\begin{picture}(\xext,\yext)(\xoff,\yoff)
 \settriparms[1`1`1;600]
  \putVtriangle(0,50)[M_1\otimes(I
  \otimes {M_2})`(M_1\otimes I)\otimes {M_2}`M_1\otimes{M_2};\alpha_{M_1, I, {M_2}}`M_1\otimes \lambda_{{M_2}}`\rho_{M_1} \otimes {M_2}]

   \put(-1400,280){$\bf MC2$}
\end{picture}
\end{center}
commute.

A {\em strict monoidal functor}
\[F:(\cM,\otimes,I,\alpha,\lambda,\rho)\lra (\cM',\otimes',I',\alpha',\lambda',\rho')\]
is a 1-cell $F:\cM\ra \cM'$ such that
\[F(I) = I', \hskip 5mm  F(M_1\otimes{M_2})= F(M_1)\otimes'F({M_2}), \]
for $M_1$, $M_2$  in $\cM$,  and
\[ F(\alpha) = \alpha'_{F\times F\times F},\hskip 5mm F(\lambda) = \lambda'_F,\hskip 5mm F(\rho) = \rho'_F.\]

A {\em monoidal transformation between two strict monoidal 1-cells}
\[\tau:F\ra F':(\cM,\otimes,I,\alpha,\lambda,\rho)\lra (\cM',\otimes',I',\alpha',\lambda',\rho')\]
is a 2-cell $\tau:F\ra F'$ in $\cK$ such that
\[ \tau_I = 1_{I'}, \hskip 5mm \tau_{M_1}\otimes\tau_{{M_2}}=\tau_{M_1\otimes {M_2}}. \]

In this way we have defined the 2-category $\Mon_{st}(\Cat)$ of monoidal categories with strict monoidal functors and monoidal natural transformations. It is an object map of a 3-functor
\[ \Mon_{st} : 2\CAT_\times \lra 2\CAT_\times, \]
whose remaining parts of definition are left for the reader.

\vskip 2mm

Given a monoidal category $(\cM,\otimes,I,\alpha,\lambda,\rho)$ in $\Cat$, we define the category of monoids. The objects are monoids, i.e., triples $(M,m:M\otimes M \ra M,i:I\ra M)$ so that the following diagram
\begin{center} \xext=2600 \yext=1250
\begin{picture}(\xext,\yext)(\xoff,\yoff)
\setsqparms[0`1`1`1;2000`600]
\putsquare(0,500)[M\otimes (M\otimes M)`M\otimes M`M\otimes M`M;`1\otimes m`m`m]

\putmorphism(0,1100)(1,0)[\phantom{M\otimes (M\otimes M)}`(M\otimes M)\otimes M`\alpha]{1000}{1}a
\putmorphism(1000,1100)(1,0)[\phantom{(M\otimes M)\otimes M}`\phantom{M\otimes M}`m\otimes 1_M]{1000}{1}a

\putmorphism(0,500)(0,1)[\phantom{M\otimes M}`I \otimes M`i \otimes 1_M]{500}{-1}l
\putmorphism(450,70)(4,1)[\phantom{I \otimes M}`\phantom{M}`]{1200}{1}r
 \put(1000,100){$\lambda_M$}

 \putmorphism(2000,1100)(1,0)[\phantom{M\otimes M}`M\otimes I`1_M \otimes i]{900}{-1}a
      \put(2830,1030){\vector(-3,-2){730}}
      \put(2430,700){$\rho_M$}

\end{picture}
\end{center}
commutes. A homomorphism of monoids $h:(M,m,i)\ra (M',m',i')$ is a morphism $h:M\ra M'$ in $\cM$ such that the diagram
\begin{center} \xext=1400 \yext=600
\begin{picture}(\xext,\yext)(\xoff,\yoff)
\setsqparms[1`1`1`1;800`500]
\putsquare(0,50)[M\otimes M`M`M'\otimes M'`M';m`h\otimes h``m']

 \settriparms[-1`0`1;500]
 \putptriangle(800,50)[\phantom{M}`I`\phantom{M'};i` `i']

 \put(710,300){$h$}

\end{picture}
\end{center}
commutes. In this way we have defined the category of monoids $\mon(\cM,\otimes,I,\alpha,\lambda,\rho)$ over monoidal category $(\cM,\otimes,I,\alpha,\lambda,\rho)$ in $\Mon_{st}(\Cat)$.  It is an object map of a 2-functor
\[ \mon : \Mon_{st}(\Cat)\lra \Cat \]
whose remaining parts of definition are again left for the reader.

\subsection{Action of monoidal category and actions along action}\label{subsec-actions-of-monoidal-objects}

A {\em monoidal action $(\cM,\otimes,I,\alpha,\lambda,\rho,\cX,\star, \psi,\bar{\psi})$ of a monoidal category $(\cM,\otimes,I,\alpha,\lambda,\rho)$  on a category} $\cX$ consists of
\begin{enumerate}
  \item a monoidal category $(\cM,\otimes,I,\alpha,\lambda,\rho)$,
  \item a category $\cX$,
  \item a functor $\star: \cM\times \cX \ra \cX$,
  \item and two natural transformations
\[ \psi: \star\circ (1_\cM\times \star )\ra \star\circ (\otimes\times 1_\cM) :\cM\times \cM\times \cX \lra cX, \]
\[ \bar{\psi}: 1_\cX\ra \star\circ \lk I, 1_X\rk:\cX\ra \cX,\]
(where $I$ denotes here the constant functor $\cM \ra\cM$ equal $I$)
\end{enumerate}
such that, for objects $M_1, M_2, M_3$ in $\cM$ and object $X$ in $\cX$,  the diagrams {\bf MA1}, {\bf MA2}, {\bf MA3} in $\cM$ below
 \begin{center} \xext=2300 \yext=1150
\begin{picture}(\xext,\yext)(\xoff,\yoff)
\put(800,1000){${M_1}\star({M_2}\star ({M_3}\star X))$}
      \put(900,950){\vector(-3,-2){380}}
      \put(1520,840){$_{\psi_{{M_1}, {M_2}, {M_3}\star X}}$} %X left up

      \put(1300,950){\vector(3,-2){380}}
      \put(100,840){$_{1_{M_1}\star\psi_{{M_2},{M_3},X}}$} %right up

\put(1450,600){$({M_1}\otimes {M_2})\star ({M_3}\star X)$} %X
\put(0,600){${M_1}\star (({M_2}\otimes {M_3})\star X)$}
      \put(500,550){\vector(1,-2){160}}
      \put(1640,330){$_{\psi_{{M_1}\otimes {M_2}, {M_3},X}}$} %left down

      \put(1700,550){\vector(-1,-2){160}}
      \put(60,330){$_{\psi_{{M_1},{M_2}\otimes {M_3}, X}}$} %right down

\put(1400,100){$(({M_1}\otimes {M_2})\otimes {M_3})\star X$} %Y
\put(50,100){$({M_1}\otimes ({M_2}\otimes {M_3}))\star X$}

  \put(1080,120){\vector(1,0){300}}
      \put(950,0){$_{\alpha_{{M_1},{M_2},{M_3}}\star 1_{X}}$} %X hor bot
\put(-700,280){$\bf MA1$}
\end{picture}
\end{center}

\begin{center} \xext=1200 \yext=750
\begin{picture}(\xext,\yext)(\xoff,\yoff)
\setsqparms[1`1`-1`1;1200`450]
 \putsquare(0,50)[{M_1}\star M_X`{M_1}\star M_X`I\star({M_1}\star X)`(I\otimes {M_1})\star X ;1_{{M_1}\star X}`\bar{\psi}_{{M_1}\star X}`\lambda_{M_1}\star 1_{X}`\psi_{I,{M_1},X}]
  \put(-1400,280){$\bf MA2$}
  \end{picture}
\end{center}
and
\begin{center} \xext=1200 \yext=600
\begin{picture}(\xext,\yext)(\xoff,\yoff)
 \put(-1400,280){$\bf MA3$}
 \setsqparms[1`1`-1`1;1200`450]
 \putsquare(0,50)[{M_1}\star X `{M_1}\star X`{M_1}\star(I\star {X})`({M_1}\otimes I)\star X;1_{{M_1}\star X}`1_{{M_1}}\star \bar{\psi}_{X}`\rho_{M_1} \star 1_{X}`\psi_{{M_1},I,X}]
 \end{picture}
\end{center}
commute.

A {\em (strict) morphism of monoidal actions}
\[ (F,G): (\cM,\otimes,I,\alpha,\lambda,\rho,\cX,\star, \psi)\ra (\cM',\otimes',I',\alpha',\lambda',\rho',\cX',\star', \psi')\]
consists of
\begin{enumerate}
  \item  a strict monoidal functor  $F: (\cM,\otimes,I,\alpha,\lambda,\rho)\lra  (\cM',\otimes',I',\alpha',\lambda',\rho')$,
  \item a functor $G: X\ra X'$
  \item such that
  \[ F(M)\star' G(X) = G(M\star X), \hskip 5mm G(\psi)=\psi'_{F\times F\times G}, \hskip 5mm G(\bar{\psi})=\bar{\psi}'_G. \]
  for objects $M$ in $\cM$ and $X$ in $\cX$.
  \end{enumerate}

A {\em transformation of strict morphisms of actions}
\[ \tau: F\ra F': (\cM,\otimes,I,\alpha,\lambda,\rho,\cX,\star, \psi)\ra (\cM',\otimes',I',\alpha',\lambda',\rho',\cX',\star', \psi')\]
is a pair of natural transformations $\tau:F\ra F'$ and $\sigma :G\ra G'$ such that
\[ \tau_{M}\star' \sigma_{X} =\tau_{M\star X}, \]
i.e.,
\begin{center} \xext=1000 \yext=700
\begin{picture}(\xext,\yext)(\xoff,\yoff)
\setsqparms[1`1`1`1;1000`600]
\putsquare(0,50)[F(M)\star' G(X)`G(M\star X)`F'(M)\star' G'(X)`G'(M\star X);=`\tau_{M}\star' \sigma_X`\tau_{M\star X}`=]
\end{picture}
\end{center}
 for objects $M$ in $\cM$ and $X$ in $\cX$.

\vskip 2mm

In this way we have defined the 2-category $\Act_{st}(\Cat)$ of actions of monoidal categories with strict morphisms and 2-cells. It is an object map of a 3-functor
\[ \Act_{st} : 2\CAT_\times \lra 2\CAT_\times, \]
whose remaining parts of definition are left for the reader.
\vskip 2mm

Given an action of monoidal category $(\cM,\otimes,I,\alpha,\lambda,\rho,\cX,\star, \psi)$ in $\Cat$, we define category $\act(\cM,\otimes,I,\alpha,\lambda,\rho,\cX,\star, \psi)$ the category of actions of monoids along the action of monoidal category $(\cM,\otimes,I,\alpha,\lambda,\rho)$ on objects of the category $\cX$. The objects are actions of monoids, i.e., 5-tuples
$(M,X,m:M\otimes M \ra M,i:I\ra M, a:M\star X\ra X )$ so that $(M,m,i)$ is a monoid and moreover the following diagram

\begin{center} \xext=2600 \yext=700
\begin{picture}(\xext,\yext)(\xoff,\yoff)
\setsqparms[0`1`1`1;2000`600]
\putsquare(0,0)[M\otimes (M\star X)`M\star X`M\star X`X;`1\otimes a`a`a]

\putmorphism(0,600)(1,0)[\phantom{M\otimes (M\star X)}`(M\otimes M)\star X`\psi]{1000}{1}a
\putmorphism(1000,600)(1,0)[\phantom{(M\otimes M)\star X}`\phantom{M\star X}`m\otimes 1_X]{1000}{1}a

 \putmorphism(2000,600)(1,0)[\phantom{M\star X}`I\star X`i\star 1_X]{900}{-1}a
      \put(2830,530){\vector(-3,-2){730}}
      \put(2520,230){$\bar{\psi}$}
\end{picture}
\end{center}
commutes. A homomorphism of actions $(h,k):(M,X,m,i,a)\ra (M',X',m',i',a')$ is a homomorphism of monoids $h:(M,m,a)\ra (M',m',i')$ and a morphism $k:X\ra X'$ in $X$ making the diagram
\begin{center} \xext=800 \yext=650
\begin{picture}(\xext,\yext)(\xoff,\yoff)
\setsqparms[1`1`1`1;800`500]
\putsquare(0,50)[M\star X`X`M'\star X'`X';a`h\star k`k`a']

\end{picture}
\end{center}
commute. In this way we have defined the category of actions of monoids $\act(\cM,\otimes,I,\alpha,\lambda,\rho,X,\star, \psi)$ over an action of a  monoidal category $(\cM,\otimes,I,\alpha,\lambda,\rho,X,\star, \psi)$ in $\Mon_{st}(\Cat)$. It is an object map of a 2-functor
\[ \act : \Act_{st}(\Cat)\lra \Cat \]
whose remaining parts of definition are again left for the reader.

\end{document}